\documentclass{article}
\usepackage{amsmath,amsthm,amsfonts}
\usepackage{graphicx,epsfig,a4,longtable}

\newcommand{\ZZ}{\mathbb{Z}}

\newcommand{\Pe}{\mathrm{Pe}}
\newcommand{\Du}{\mathrm{Du}}

\newcommand{\Aut}{\mathrm{Aut}}
\newcommand{\RR}{\mathit{R}}
\newcommand{\TT}{\mathit{T}}
\newcommand{\LL}{\mathit{L}}

\newcommand{\Mon}{\mathrm{Mon}}
\newcommand{\Orient}{\mathrm{Or}}

\newcommand{\id}{\mathrm{\underline{id}}}
\newcommand{\idd}{\mathrm{\underline{\underline{id}}}}
\newcommand{\Id}{\mathrm{Id}}
\newcommand{\Symr}{\mathrm{Sym_R}}
\newcommand{\Syml}{\mathrm{Sym_L}}
\newcommand{\Flags}{\mathrm{Flags}}
\newcommand{\Orb}{\mathrm{Orb}}
\newcommand{\Lifts}{\mathrm{Lifts}}
\newcommand{\Co}{\mathrm{Co}}
\newcommand{\DM}{\mathrm{DM}}
\newcommand{\lcm}{\mathrm{lcm}}
\newcommand{\Dih}{\mathrm{D}}
\newcommand{\MN}{\mathrm{MN}}
\newcommand{\dd}{\mathbf{du}}
\newcommand{\pp}{\mathbf{pe}}
\newcommand{\Type}{\mathbf{Type}}
\newcommand{\Ga}{G^1}
\newcommand{\Gb}{G^2}
\newcommand{\Gi}{G^i}
\newcommand{\Ha}{H^1}
\newcommand{\Hb}{H^2}
\newcommand{\Hi}{H^i}
\newcommand{\MonQ}{\triangle}

\theoremstyle{definition}

\theoremstyle{plain}

\newtheorem{lemma}{Lemma}[section]
\newtheorem{theorem}[lemma]{Theorem}
\newtheorem{proposition}[lemma]{Proposition}
\newtheorem{corollary}[lemma]{Corollary}

\newtheorem{example}[lemma]{Example}

\usepackage{amsmath,amssymb}
\renewcommand{\title}[2][]{
  \begin{center}
    \begin{large}
      \begin{bf}
        #2
      \end{bf}
    \end{large}
  \end{center}}

\renewcommand{\author}[2][]{
  \begin{center}
    \begin{sc}
      #2
    \end{sc}
  \end{center}}

\newcommand{\address}[1]{
  \begin{center}
        #1
  \end{center}}

\newcommand{\email}[1]{
  \begin{center}
    \begin{tt}
      #1
    \end{tt}
  \end{center}
}

\renewcommand{\maketitle}{}

\begin{document}

\title{Parallel-product decomposition of edge-transitive maps}

\author[Orbanic]{Alen Orbani\'c}
\address{
University of Ljubljana \\
Faculty of Mathematics and Physics\\
Department of Mathematics \\
Jadranska 19 \\
1000 Ljubljana, SLOVENIA
}

\email{Alen.Orbanic@fmf.uni-lj.si}

\date{\today}

\maketitle

\begin{abstract}
The parallel product of two rooted maps was introduced by S. E. Wilson in 1994. 
The main question of this paper is whether for a given reflexible map $M$ one can decompose
the map into a parallel product of two reflexible maps. This can be achieved if and only if the 
monodromy (or the automorphism) group of the map has at least two minimal normal subgroups.
All reflexible maps 
up to 100 edges, which are not parallel-product decomposable, are calculated and presented. For this purpose, all
degenerate and slightly-degenerate reflexible maps are classified.

Three different quotients of rooted maps are considered in the paper and 
a characterizaton of morphisms of rooted maps similar to the first isomorphism theorem for groups is presented.
The monodromy quotient of a map is introduced, having the property 
that all the automorphisms project. 

A theory of edge-transitive maps on non-orientable surfaces is developed. A concept of reduced regularty in the 
manner of Breda d'Azevedo is applied on edge-transitive maps. Using that,
the concept of parallel-product decomposability is extended to edge-transitive maps, where a  
characterization in terms of minimal normal subgroups of the automorphism group is obtained.
Additionally, using Petrie triality and 
the parallel-product decomposition, a new organization of edge-transitive maps is presented, providing a basis for future
censuses.
\end{abstract}


\section{Introduction}

{\bf History and motivation.} The history of edge-transitive maps, which also include regular (reflexible) and orientably regular maps, starts with ancient Greeks (the platonic solids, also some of the archimedian solids). 
In the 17th century, Kepler \cite{kepler} worked on stellated polyhedra where some non-planar regular maps occured.
In the 19th century, Heffter \cite{heffter} considered orientably regular embeddings of complete graphs, while 
Klein \cite{klein} and Dyck \cite{dyck} 
constructed some cubic regular maps on the surface of orientable genus 3, in the context of automorphic functions. In the beginning of the 20th century, regular maps were
first used as geometrical representations of groups (Burnside \cite{burnside}). 
More systematic study of regular maps continued with 
Brahana \cite{brahana} and Coxeter and Moser \cite{coxeterMoser}, where regular maps were treated as geometrical, 
combinatorial and group theoretical objects. The basis for the modern treatment of general maps was set 
by Jones and Singerman \cite{jonesSingerman} for orientable surfaces  and 
by Bryant and Singerman \cite{bryantSingerman} for non-orientable ones. The classic reference for maps became the book
by Gross and Tucker \cite{grossTucker}. In the last decade, research on maps of high symmetry has mainly
focused on regular (and orientably regular) maps and Cayley maps. The recent paper by
Richter, \v Sir\'a\v n, Jajcay, Tucker and Watkins \cite{cayleyMaps} provides a nice survey for Cayley maps.

For edge-transitive maps, Graver and Watkins \cite{memoirs} give 
the fundamental classification into 14 types according to the possesion of some types of automorphisms. 
The existence of all of the types on infinitely many orientable surfaces was shown in the 
important work by \v Sir\'a\v n, Tucker and Watkins \cite{realizing}.

The central problem of edge-transitive maps is construction and  classification. 
The most common constructions of edge-transitive maps arise either from constructions of 
finite groups admitting one of 14 types of presentations \cite{realizing}  
or as covers of smaller maps.
Three natural approaches are used in the classification of edge-transitive maps, namely 
by the number of edges, by the  
 underlying surface and  by the underlying graph. The results of those classifications are several censuses
\cite{lowx, families, census}.

It is known that 
all compact closed surfaces, other than the sphere, torus, projective plane and Klein bottle, necessarily contain a finite number 
of edge-transitive maps. The upper bound depends on the surface and it is easily obtained from Euler's formula or the Riemann-Hurwitz equation.
All edge-transitive maps on the torus were classified by \v Sir\'a\v n, Tucker and Watkins \cite{realizing},
the classification for the sphere was done by Gr\" unbaum and Shephard \cite{grunbaumShephard}, while 
a part of the classification for the Klein bottle was done by Poto\v cnik and Wilson \cite{potocnikWilson}.

Before the age of fast computers,
many authors (Brahana \cite{brahana}, Coxeter and Moser \cite{coxeterMoser},
Sherk \cite{sherk}, Garbe \cite{garbe}, Bergau and Garbe \cite{bergauGarbe}) worked on the classification of regular and
orientably regular maps and managed to
classify all regular and orientably regular maps on surfaces of orientabe genus up to  7 and non-orientable genus up to 8. 
In the 1970s, Wilson in his Ph.D. thesis \cite{wilsonThesis} calculated most reflexible and chiral 
maps up to 100 edges \cite{census} using a computer and running his {\em Riemann surface algorithm} \cite{riemannAlg}.  
The recent breakthrough in this field is due to Conder and Dobcs\' anyi \cite{smallgenus}, 
who calculated all orientably regular maps on surfaces from genera 3 up to 
15 and all the non-orientable reflexible maps on surfaces from non-orientable genera 2 up to 30 
(Conder\&Dobcs\' anyi's census \cite{lowx}).

Since Wilson's and Conder\&Dobcs\' anyi's censuses present different information, a census was needed that would contain the
information from both of them. Since chiral maps are not closed under the Petrie dual, the natural extension
of the censuses, as observed by T. Pisanski, seemed to be edge-transitive maps. 
Such an extended census was the motivation for this paper. Since the census is so large, 
the author sought a shorter description
of maps in terms of some kind of ''primitive'' maps from which all other 
maps can be obtained using some
set of operations. The algorithms for performing the operations needed to be of 
relatively low time complexity so the computations
of ''non-primitive'' maps remain simple.  
It turned out that the appropriate operation is the parallel product introduced by Wilson \cite{wilsonParallel}. 

{\bf Overview of main results.} Let us focus on some special class of edge-transitive maps -- reflexible maps. The question is
which reflexible maps are {\em parallel-product decomposable}, that is, a parallel product of two reflexible maps. The maps that are not parallel-product decomposable 
are of special interest as  basic building blocks.

Usually, a map is represented by a set of flags and by three involutions, two of which commute, treated as permutations of the flags and intuitively giving instructions for gluing the flags together to form a surface.
The group generated by these three involutions acts transitively on the set of flags and is called the {\em monodromy group} of the map. The automorphism group of
a map is the group of permutations of the flags respecting the action of the monodromy group.

The main results of this article are the following group theoretical characterizations of parallel-product decomposability:\\[2mm]
{\bf Theorem \ref{reflexibleDecomposition}.}
A reflexible map is parallel-product decomposable if and only if the monodromy group (and therefore also the  
automorphism group) contains at least two non-trivial minimal normal subgroups. \\[2mm]
{\bf Theorem \ref{theFinal}.}
An edge-transitive map $M$ is parallel-product decomposable if and only if $\Aut(M)$ contains at least two minimal normal
subgroups. \\[2mm]
\indent These two theorems are consequences of the main result of the paper:\\[2mm]
{\bf Theorem \ref{decomposition}.}(Decomposition theorem) A map $M = (f, G, Z, \id)$  is parallel-product decomposable
if and only if there are at least two normal non-transitive subgroups $\Ha, \Hb \vartriangleleft G$, such that 
$\Ha \cap \Hb = \{1\}$ and $G_{\id}\Ha \cap G_{\id}\Hb = G_{\id}$.\\[2mm]

{\bf Paper layout.} The sections of this paper are organized as follows. Section 2 establishes the algebraic machinery necessary to discuss rooted maps
on all surfaces in a manner similar to the article about Cayley maps
by Richter, \v Sir\'a\v n, Jajcay, Tucker and Watkins \cite{cayleyMaps}. The algebraic machinery includes also
rooted map morphisms and vertex-face-Petrie circuits triality.

Section 3 contains results about quotients of maps. Only one type of quotient of maps
has appeared in the literature, namely a quotient here called an {\em automorphism quotient} defined by Malni\v c, Nedela and 
\v Skoviera 
\cite{malnicSkovieraNedela}, or regular covering map in \cite{grossTucker}. Here we introduce two completely new 
quotients: a $K$-quotient and a monodromy quotient. A complete characterization of map morphisms in terms of quotients 
in the manner of the first isomorphism theorem on groups is given by Theorem \ref{charQuo}. 
The {\em monodromy quotient}, which has the important property that all automorphisms project,  is introduced at 
the end of the section.

Section 4 describes some interesting properties of the parallel product. 
After Wilson\cite{wilsonParallel} introduced the parallel product, only a few authors considered it as an important 
operation on maps (see \cite{joinIntersection} for hypermaps). The most interesting result of this 
section is the construction of the 
smallest unique reflexible cover above any map. Similarly, the unique totally symmetric cover of a reflexible map 
is also constructed.
Lifts of automorphisms in a parallel product of maps are studied.

Section 5 is devoted to the proof of Theorem \ref{decomposition}.

Section 6 classifies all degenerate and slightly degenerate reflexible maps. 
These are basically the maps containing vertices of valence less than 3 or some kind of a 
degeneracy of the edges, such as loops or semi-edges. These degeneracies arise naturally in quotients and in the {\em triality} of duals and Petrie duals.

In Section 7 all 
parallel-product indecomposable reflexible maps up to 100 edges are presented. The most important theorem of 
this section is the decomposition theorem for 
reflexible maps, namely Theorem \ref{reflexibleDecomposition}.

Section 8 extends the theory of edge-transitive maps introduced in \cite{memoirs, realizing} 
to non-orientable surfaces and 
presents the organization of a census of edge-transitive maps using triality and the parallel product.
By triality, the 14 types are reduced to the 6 basic types needed for the reconstruction of all 
edge-transitive maps. 
For these
types, partial finite presentations of the corresponding automorphism groups are given, as well as a method to uniquely 
reconstruct the corresponding maps. Using the concept of {\em reduced regularity} 
introduced by Breda d'Azevedo \cite{reducedBreda}, 
a theory of presentations for edge-transitive maps, is developed. 
Here we change the presentation of a map, and thus the monodromy group,
so that both the new monodromy group and the automorphism group become regular. This approach enables us to
use a characterization of parallel-product decomposability for edge-transitive maps, namely Theorem \ref{theFinal},  
and forms a basis for  future work.

\section{Definitions}

In the present work we will denote by $\Symr(S)$ a symmetric group on $|S|$ elements, i.e. the set 
of all the permutations
on the elements of $S$, such that the composition of the permutations is done from the left to the right. Also, $\Symr(S)$ naturally
acts on the set $S$ by the right action. Similarly, we will denote by $\Syml(S)$ 
the set of all the permutations of the elements 
of the set $S$, but here we have the composition from the right to the left, as functions are usually composed. 
Also, $\Syml(S)$
naturally acts on the set $S$ by the left action.

Let $F = \langle t, l, r\ \vert\ t^2 = l^2 = r^2 = (tl)^2 = 1\rangle$. 
A {\em (finite) rooted map} $M$ 
is a quadruple $M = (f, G, Z, \id) = (f_M, \Mon(M), \Flags(M), \id_M)$, where $Z$ is 
a finite set of {\em flags}, $G \leq \Symr(Z)$ acts transitively and faithfully  
on $Z$, $f: F \to G$ is a group epimorphism and $\id_M \in Z$ a {\em root flag}. 
Define $T := f(t), L:= f(l), R := f(r)$. The group $G = \Mon(M)$ is  referred to as  the
{\em monodromy group} of the rooted map $M$. 
We will often denote an empty word from $F$ or $\Mon(M)$ by $\epsilon$, but sometimes also by 1. 
An identity mapping on a set $S$ is often denoted by $\Id$.

Note that a monodromy group as an algebraic object is not just a group, but a group together with the labelled 
generators ($\TT$, $\LL$ and $\RR$). For two groups $G$ and $K$ generated by $k$ generators 
labelled with labels $a_1, \ldots, a_k$, we will say that $G$ and $K$ are {\em congruent} if there 
exists an isomorphism of $G$ and $K$, which respects the labelling, and therefore maps a generator of $G$
labelled by $a_i$ to the generator of $K$ also labelled by $a_i$, for $i = 1, \ldots, k$.
We will denote the congruence by $G = K$ and the corresponding isomorphism will be called the 
{\em congruence isomorphism}.
When the groups we are working with are monodromy groups, the labels to be considered are $\TT$, $\LL$ and $\RR$.
 
If we ignore the choice of a root flag, we obtain maps named {\em holey maps} used in \cite{consForbLift}.  
The word ''map'' will refer to  the word ''rooted map'' in this work. In general, a right action
of a group $G$ on a set $Z$ will be denoted by $(Z, G)$ and a  left action by $(G, Z)$. We will denote a stabilizer of an
element $z \in Z$ by $G_z$.

The {\em flag graph} $\Co(M)$ is the trivalent 
multigraph with a vertex set 
$\Flags(M)$, where each $x\in \Flags(M)$ is connected with flags $x\cdot \TT, x\cdot \LL, x\cdot\RR$.
If for $x \in Z$, $W \in \{\TT, \LL, \RR\}$, it follows $x\cdot W = x$, then a semi-edge
emanates from $x$. The involutions $\TT$, $\LL$ and $\RR$ 
naturally induce a 3-coloring of edges and semi-edges of $\Co(M)$.

Let $M$ and $N$ be two rooted maps. A morphism of the maps is a pair $(\phi,\psi)$, where 
$\psi: \Mon(M) \to \Mon(N)$ is an epimorphism of the groups, such that $\psi\circ f_M = f_N$
and $\phi:\Flags(M) \to \Flags(N)$ is an
onto mapping, where $\phi(\id_M) = \id_N$ and $\phi(z\cdot g) = \phi(z)\cdot \psi(g)$ 
for every $z \in \Flags(M)$ and $g \in \Mon(M)$. 
A morphism of rooted maps is also called a {\em covering projection}. In this case the map $M$ is called a {\em cover} of 
the map $N$. 
Note that the notion of covering projection corresponds to the notion of covering projection of flag graphs described
in \cite{liftMalnic}. Since such a projection can take an edge to a semi-edge, this kind of a 
projection is not a covering projection in the sense of topology, namely a local homeomorphism, 
but is more like the projection associated with an orbifold.

If both $\phi$ and $\psi$ are one-to-one then the pair $(\phi, \psi)$ is an {\em isomorphism} of the rooted maps. 
If we omit the condition $\phi(\id_M) = \id_N$ then we get a {\em generalized isomorphism} of rooted maps.
Note that this is an isomorphism of holey maps.

An {\em automorphism} of a rooted map $M$ is a generalized isomorphism $(\phi, \Id): M \to M$, 
where $\Id$ denotes the identity 
mapping of $\Mon(M)$.
The group of all automorphisms is denoted by $\Aut(M)$. 
Since for any $W \in \Mon(M)$, $\alpha \in \Aut(M)$, $x \in \Flags(M)$, it follows $\alpha(x\cdot W) = \alpha(x)\cdot W$,
each automorphism is already defined by a mapping of a single flag. Thus $\Aut(M)$ acts semi-regularly on 
$\Flags(M)$. A map $M$ is {\em reflexible} if and only if $\Aut(M)$ is regular. It is well known that this is true for maps 
if and only if $\Aut(M) \simeq \Mon(M)$ and $\Mon(M)$ is also regular. In a slightly general form this will be also 
proved in Proposition \ref{monAutReg}.
Given $W \in F$, we say that a rooted map $M$ {\em contains} the automorphism $\alpha_W$, if there is  
an automorphism of $M$
taking the flag $\id$ to the flag $\id \cdot f(W)$. 
If two maps 
contain $\alpha_W, W \in F$, we will say that the maps have $\alpha_W$ in {\em common}.

The edges $E(M)$ of a map $M$ are the orbits of $\langle \TT, \LL \rangle$, where $\langle \TT, \LL\rangle$ denotes the subgroup
of $\Mon(M)$ generated by $\TT$ and $\LL$. The vertices $V(M)$ are the 
orbits of $\langle \TT, \RR\rangle$, the faces $F(M)$ are the orbits of $\langle \LL, \RR\rangle$ and the 
Petrie circuits $P(M)$ are the orbits of $\langle \TT\LL, \RR\rangle$. Let 
$\Orient(M) := \langle \RR\TT, \RR\LL \rangle$ denote the image in $\Mon(M)$ 
of the index two subroup of $F$ consisting of even length words. 
It is easy to see that the number of orbits of the action $(Z, \Orient(M))$ 
is 1 or 2. 
It is known, that in the case when $\TT,\LL,\RR$ are not contained in any stabilizer  
of any flag (i.e. they are fixed-point-free),
the map
combinatorially represents a map on a compact closed surface. 
In in this case, we say that the map is {\em orientable} if $\Orient(M)$ has 2 orbits, and 
{\em non-orientable} otherwise. 
If $\TT\LL$ has a fixed point then the map has a semi-edge. 
Note that if $\Co(M)$ does not have semi-edges, then orientability coincides with
with the graph $\Co(M)$ being bipartite.

The {\em parallel product} of two maps $M = (f_1, G_1, Z_1, \id_1)$ and $N = (f_2, G_2, Z_2, \id_2)$ 
is defined as $M\parallel N := ((f_1,f_2), G, Z, (\id_1, \id_2))$ where 
$G := (f_1,f_2)(F) \leq G_1 \times G_2$ and $Z := \Orb_{(\id_1, \id_2)}^{G}(Z_1 \times Z_2)$ an
orbit 
of the action $(Z_1\times Z_2, G)$ containing $(\id_1, \id_2)$. 
Thus the monodromy group  $\Mon(M \parallel N)$ is a subgroup of $\Mon(M) \times \Mon(N)$ generated by 
$(\TT_M, \TT_N)$, $(\LL_M, \LL_N)$ and  $(\RR_M, \RR_N)$ and the flags are the subset of $\Flags(M) \times \Flags(N)$
that is an orbit of $\Mon(M\parallel N)$ containing the new root $(\id_1, \id_2)$. 
We will often denote a new root by $\id_{1,2}$ or $\id_{M,N}$.
It is easy to see that the action of the new monodromy group is faithful. 
A pair $(f_1, f_2)$ will be often denoted by $f_{1,2}$ or $f_{M,N}$ and similarly the set $Z_1 \times Z_2$ by
$Z_{1,2}$. Note that the parallel product is associative and commutative (up to isomorphism of the obtained maps).
This was already noted by Wilson \cite{wilsonParallel}, where the product was introduced.
A parallel product is said to be {\em non-trivial} if and only if it is not isomorphic to one of the factors. 

Let $\dd, \pp$ be automorphisms of $F$ defined by $\dd: t\mapsto l, l\mapsto t, r \mapsto r$ and 
$\pp: t \mapsto t, l\mapsto lt, r\mapsto r$. Then the {\em dual} of a map $M = (f, G, Z, \id)$
is defined as $\Du(M) := (f\circ \dd, G, Z, \id)$. The {\em Petrie dual} is defined as
$\Pe(M) := (f\circ \pp, G, Z, \id)$. 
It should be noted that given a map $M$, both $\Du(M)$ and $\Pe(M)$ have the same edges as $M$, but $\Du(M)$
interchanges faces and vertices leaving Petrie circuits the same, while
$\Pe(M)$ interchanges faces and Petrie circuits leaving vertices the same.
Since $\langle \dd, \pp\rangle \simeq S_3$, as a subgroup of $\Aut(F)$,
at most $6$ non-isomorphic
maps can be produced applying these two operations. 
Since all the maps obtained using the operations $\Du$ and $\Pe$ have 
the same automorphism group (only the roles of automorphisms are changed), we will often analyze only one 
representative of the class. 
The symmetry provided by $\langle \dd, \pp\rangle$ will be called a {\em triality} and a class of maps obtained
from a single map by applying the operations will be called a {\em triality class}.

Let $p: \widetilde{X} \to X$ be a morphism of maps. Let 
$\widetilde{f}\in \Aut(\widetilde{X})$. If there exists
$f \in \Aut(X)$, such that $p\circ \widetilde{f} = f\circ p$, then we say
that $\widetilde{f}$ {\em projects} (along $p$). 
On the other hand, if there is $f\in \Aut(X)$ and there exists 
$\widetilde{f}\in\Aut(\widetilde{X})$, such that $p\circ \widetilde{f} = f\circ p$,
we say that $f$ {\em lifts} (with $p$) and 
$\widetilde{f}\in\Lifts_p(f)$ is one of its lifts. Note that for $W \in F$, if 
$\widetilde{\alpha}_W \in \Aut(\widetilde{M})$ projects, it projects to 
$\alpha_W$.

If a root flag $\id$ of a map $M$ is changed to the flag $\id\cdot W$, $W \in \Mon(M)$,
a {\em re-rooted} map is obtained. If $W \in \{\epsilon, \TT, \LL, \TT\LL\}$, the 
obtained re-rooted map is said to be {\em simply re-rooted}.
In general, re-rooted maps are not isomorphic as rooted maps, although they are isomorphic as
holey maps.

\section{Quotients of maps}

In this section, for an arbitrary map $M$, three different quotients are introduced, namely a $K$-quotient, 
for some subgroup $K \leq \Mon(M)$,  
a monodromy quotient and an automorphism quotient. It is shown that 
any image of a map by a map morphism is isomorphic to some $K$-quotient.  

The topics discussed in this generalize the work of Malni\v c, Nedela and \v Skoviera \cite{malnicSkovieraNedela}. They mainly worked with quotients of a regular map obtained through
subgroups of the automorphism group, while in this paper we mainly work with quotients obtained through subgroups of the monodromy group. 

Let $(Z, G)$ be a transitive action. Then all the stabilizers are conjugate and their intersection is
a normal subgroup $H \vartriangleleft G$. Let $\chi: G \to S_{|Z|}$ be the homomorphism of groups mapping 
$g \in G$ to the corresponding permutation that acts on the elements of $Z$ in the same manner as
$g$. It is easy to see that $\ker \chi = H$. Since $G/H$ is isomorphic to $\chi(G)$, the isomorphism induces
an action $(Z, G/H)$ defined by $z\cdot Hg = z\cdot g$, for any $z \in Z$ and $g \in G$, where $Hg \in G/H$.  
Since the action $(Z, \chi(G))$ is faithful, the action $(Z, G/H)$ is also faithful.
In this case $H$ is called the {\em kernel of the action $(Z, G)$}. 

The following proposition defines a way of obtaining the first kind of quotient of a map.

\begin{proposition}
\label{KQuoDef}
Let $M = (f, G, Z, \id)$ and let  $K \leq G$ be a subgroup, such that $G_{\id} \leq K$. 
Let $H$ be the kernel of the action $(G/K, G)$ and $q: G \to G/H$ be the natural epimorphism. 
Then $N = (q \circ f, G/H, G/K, K)$ is a map and there exists a map morphism $(p,q): M \to N$.
\end{proposition}

\begin{proof}
Since the action 
$(G/K, G/H)$ is transitive and faithful, $N$ is a map.
Note that since $H \vartriangleleft K$, the action $(G/K, G/H)$ is naturally defined by
$Ka\cdot Hb = KaHb  = K(aHa^{-1})ab = Kab$. 
Define $p: Z \to G/K$ by $p(\id\cdot g) = Kg$ for any $g \in G$. 
Let $x \in Z$ and $g, h \in G$, such that $x = \id\cdot g = \id\cdot h$. 
Then $gh^{-1} \in G_{\id}\leq K$ and $p(x)$ is well defined.

Let $z \in Z$ and $g \in G$ be arbitrary and $h \in G$, such that $z = \id\cdot h$.
Then $p(z \cdot g) = p(\id \cdot hg) = Khg$. Also, $p(z)\cdot q(g) =$ $p(\id\cdot h)\cdot q(g) =$ 
$Kh\cdot Hg = Khg$. As $p(\id) = K$ and $p$ is obviously onto, 
$(p, q)$ is a map morphism.
\end{proof}


Any map $N$ obtained from $M$ in the way shown in Proposition \ref{KQuoDef} is called a {\em $K$-quotient} 
and denoted by $M/K$. The following theorem characterizes all the images of morphisms of a given map.

\begin{theorem}
\label{charQuo}
Let $M = (f, G, Z, \id)$, $N = (f_N, G_N, Z_N, \id_N)$ and let
$(\phi, \psi): M \to N$ be a map morphism. Then $N$ is isomorphic to $M/K$, 
where $K = \psi^{-1}\left((G_N)_{\id_N}\right)$ and 
$(G_N)_{\id_N}$ denotes the stabilizer of $\id_N \in Z_N$ of the action $(Z_N, G_N)$.
In particular, every image $N$ of any map morphism from $M$ is isomorphic to some   
$K$-quotient for $G_{\id}\leq K \leq G$.
\end{theorem} 

\begin{proof}

By the definition of a map morphism, $f_N = \psi\circ f$. Since $\psi$ is onto, $G_N$ is isomorphic
to $G/H$, where $H = \ker \psi$. Let $s: G_N \to G/H$ be that isomorphism and $q: G \to G/H$ a natural
epimorphism. Then $q = s\circ \psi$. 

Let $K := \psi^{-1}((G_N)_{\id_N})$. Thus $G_{\id} \leq K \leq G$. 
The stabilizer of the coset $K$ of the action $(G/K, G)$ is exactly $K$. The kernel of the action is thus the
intersection of all the conjugates of $K$:

\begin{align*}
\bigcap_{a\in G}a^{-1}Ka &= \bigcap_{a\in G} a^{-1}\psi^{-1}((G_N)_{\id_N})a = 
\bigcap_{a\in G_N} \psi^{-1}(a^{-1}(G_N)_{\id_N}a)\\
&= \psi^{-1}\left( \bigcap_{a\in G_N}a^{-1}(G_N)_{\id_N}a\right) = \psi^{-1}(\{1\}) = H.
\end{align*}
Note that the calculation above is true because $\psi$ is onto and 
the action $(Z_N, G_N)$ is faithful and thus the kernel of the action 
equals $\bigcap_{a\in G_N}a^{-1}(G_N)_{\id_N}a = \{1\}$. 
Thus $M/K = (\psi \circ f, G/H, G/K, K)$.

For $z \in Z_N$ and $g \in G_N$ define a mapping $r: Z_N \to G/K$ by 
 $r: \id_N\cdot g \mapsto Ku$, where 
$u$ is any element from $\psi^{-1}(g)$. If $u'\in \psi^{-1}(g)$ is any other such element, then 
$u^{-1}u' \in \ker \psi \leq K$ and thus the definition is independent of the choice of $u$. If 
$\id_N\cdot g = \id_N \cdot h$, then $gh^{-1} \in G_{\id_N}$. 
Let $u \in \psi^{-1}(g)$ and $v \in \psi^{-1}(h)$. Then $uv^{-1} \in \psi^{-1}(gh^{-1}) \leq K$ and
$Ku = Kv$. Hence the mapping $r$ is well defined.

Let $z \in Z_N$, $g\in G_N$ be arbitrary and let $h \in G_N$, such that $\id_N\cdot h = z$. Let $u \in \psi^{-1}(g)$ and 
$v \in \psi^{-1}(h)$. Then $s(g) = Hu$ and $s(h) = Hv$.
Hence, $r(z\cdot g) = r(\id_N\cdot hg) = Kvu$. On the other hand, $r(z)\cdot s(g) =$ 
$r(\id_N\cdot h)\cdot s(g) =$ $Kv\cdot Hu = Kvu$.  

For $x, y \in Z_N$, $g, h \in G/H$ and $u \in \psi^{-1}(g)$, $v \in \psi^{-1}(h)$, let $x = \id\cdot g$ and 
$y = \id\cdot  h$.
Then $r(x) = r(y)$ means $uv^{-1}\in K$ implying that $gh^{-1} \in (G_N)_{\id_N}$ and $x = y$. 
Therefore $r$ is one-to-one and since it is
always onto, it is a bijection.
As $r(\id_N) = K$, the
mapping $(r, s): N \to M/K$ is a map isomorphism. 
\end{proof}

\begin{corollary}
\label{isomorph}
A map $M = (f, G, Z, \id)$ is isomorphic to its $G_{\id}$-quotient $(f, G,$ $ G/G_{\id}, G_{\id})$.
\end{corollary}
\begin{proof}
Take $(\phi, \psi) = (\Id, \Id): M \to M$ and apply Theorem \ref{charQuo}. 
\end{proof}

Another type of a quotient can be obtained in the following way. 

\begin{proposition}
\label{monQuoDef}
Let $M = (f, G, Z, \id)$ and $H \vartriangleleft G$ be a normal subgroup. Let $q: G \to G/H$ be the 
natural epimorphism. Let $Z/H$ denote the set of orbits of the action $(Z, H)$ and let $p: Z \to Z/H$ be defined 
as $p: z \mapsto [z]$, where $[z]$ denotes the orbit containing the element $z \in Z$.
Then $N = (q\circ f, G/H, Z/H, p(\id))$ is a map isomorphic to the $K$-quotient $M/G_{\id}H$. 
\end{proposition}

\begin{proof}
Note that for a word $W \in G/H$ there exists a word $w \in G$, such that $W = q(w)$. 
For an orbit $[x] \in Z/H$ we define an operation as $[x]\cdot W := [x\cdot w]$. 
For any other word $v \in \Mon(M)$, such that
$q(v) = W$, it follows $Hv = Hw$ or $vw^{-1} = h \in H$. Since $H$ is normal, 
$x\cdot v = x\cdot hw = x\cdot w(w^{-1}hw) = x\cdot wh'$,
where $h' \in H$. Therefore $[x\cdot w] = [x\cdot v]$ and the operation above is well defined. 
It is easy to verify that 
the operation indeed meets the conditions to be a right action. 
The action $(Z/H, G/H)$ is obviously transitive. 

Let $g \in G/H$, such that $g$ stabilizes $[\id]$. There exists some 
$u \in G$, such that $q(u) = g$. 
Since $[\id]\cdot g = [\id\cdot u] = [\id]$, there exists some $h \in H$, such that 
$\id\cdot u = \id \cdot h$. Therefore, $uh^{-1} \in G_{\id}$ and $q(uh^{-1}) \in HG_{\id} = G_{\id}H$. 
Note that since $H$ is normal, $G_{\id}H$ is a subgroup of $G$. 

For $g \in G_{\id}H$, it follows $g = sh$,
for some $s \in G_{\id}$ and $h\in H$. 
Then $[\id]\cdot q(g) = [\id\cdot sh] = [\id]$. 
Hence, the stabilizer of $[\id]$ is exactly $q(G_{\id}H)$. 

Since the action $(Z, G)$ is faithful, the kernel of the action $\bigcap_{x\in G}x^{-1}G_{\id}x$ is trivial.
Therefore 
\begin{align*}
\bigcap_{x\in G} x^{-1}G_{\id}Hx = \left(\bigcap_{x\in G}x^{-1}G_{\id}x\right)H = H,
\end{align*}
since $H$ is normal.
Thus $M/G_{\id}H = (q \circ f, G/H, G/G_{\id}H, G_{\id}H)$. 
For the kernel of the action $(Z/H,$ $G/H)$ it follows
\begin{align*}
\bigcap_{a\in G/H}a^{-1}q(G_{\id}H)a &= \bigcap_{x\in G}q(x)^{-1}q(G_{\id}H)q(x) = q\left(\bigcap_{x\in G} x^{-1}G_{\id}Hx\right) = q(H) = 1,\\
\end{align*}
since $q$ is onto.
Thus, $(Z/H, G/H)$ is faithful and $N$ is a map.

Let $r: Z/H \to G/G_{\id}H$ be a mapping defined by $r: [\id\cdot u] \mapsto G_{\id}Hu$. 
Since $[\id\cdot u] = [\id \cdot v]$, for some $u,v \in G$, if and only if $uv^{-1} \in G_{\id}H$, the mapping $r$ 
is well defined and one-to-one.
Obviously, it is also onto. Let $Hu \in G/H$, for some $u \in G$, and let $v \in G$. Then
$r([\id\cdot v])\cdot Hg = G_{\id}HvHg = G_{\id}Hvg$. On the other hand, 
$r([\id\cdot v]\cdot Hg) = r([\id\cdot vg]) = G_{\id}Hvg$.
As $r([\id]) = G_{\id}H$, it follows that $(r, \Id): N \to M/G_{\id}H$ is a map isomorphism.
\end{proof}

The quotient defined in Proposition \ref{monQuoDef} is called the {\em monodromy quotient induced by $H$}. 
We denote the monodromy quotient of a map $M$ induced by a normal subgroup $H \vartriangleleft \Mon(M)$ by
$M\MonQ H$.
The corresponding projection
is called the {\em monodromy quotient projection}.

The following proposition presents one of the most important properties of the monodromy quotient.

\begin{proposition}
\label{monProject}
Let $\widetilde{X} = (f, G, Z, \id)$ be a map, $H\vartriangleleft G$ a normal subgroup and $X = \widetilde{X}\MonQ H$ be the
monodromy 
quotient induced by $H$. Let $(p,q)$ be the monodromy quotient projection  
and $\widetilde{a} \in \Aut(\widetilde{X})$. Then $\widetilde{a}$ projects. 
In particular, if for $W \in F$, the map $\widetilde{X}$ contains $\alpha_W$ then the map
$X$ also contains $\alpha_W$.
\end{proposition}
\begin{proof}
Define $a([x]) := [\widetilde{a}(x)]$.
Let $y \in [x]$. Then there exists $h\in H$, such that $y = x\cdot h$ and 
$\widetilde{a}(y) = \widetilde{a}(x\cdot h) = $ $\widetilde{a}(x)\cdot h \in [\widetilde{a}(x)]$.
Thus the mapping $a$ is well defined. For $W \in q(G)$, there exists $g \in G$, such that $q(g) = W$. Then
$a([x]\cdot W) = a([x]\cdot q(g)) = a([x\cdot g]) = [\widetilde{a}(x\cdot g)] = [\widetilde{a}(x)\cdot g] = $ 
$[\widetilde{a}(x)]\cdot q(g) = a([x])\cdot W$. 

If $a([x]) = a([y])$, then $[\widetilde{a}(x)] =[\widetilde{a}(y)]$ and $\widetilde{a}(x) = \widetilde{a}(y)\cdot h = \widetilde{a}(y\cdot h)$ for some 
$h \in H$. Thus $x = y\cdot h$ and $[x] = [y]$, implying that $a$ is one-to-one. 
Obviously, it is also onto and thus $a \in \Aut(X)$. 

If for $W \in F$, $\widetilde{a} = \alpha_W \in \Aut(\widetilde{X})$ then 
$\widetilde{a}(\id_{\widetilde{X}}) = \id_{\widetilde{X}}\cdot f(W)$. 
Therefore $a(\id_X) = a([\id_{\widetilde{X}}]) = [\widetilde{a}(\id_{\widetilde{X}})] =$
$[\id_{\widetilde{X}}\cdot f(W)] = $
$[\id_{\widetilde{X}}]\cdot q(f(W)) = \id_{X}\cdot (q\circ f)(W)$ meaning that $a = \alpha_W \in \Aut(X)$.
\end{proof}

An interesting observation made by Tucker \cite{tuckerPrivate} is that any map morphism  $(\phi, \psi): M \to N$ 
factors through a monodromy quotient of $M = (f, G, Z, \id)$ obtained using $H = \ker \psi$. Let
$(p,\psi): M \to M\MonQ H$ be the monodromy quotient projection. Then 
$(\phi, \psi) = (r, \Id)\circ (p, \phi)$, where $r$ is uniquely defined by $\phi = r\circ p$, since 
$\phi$ and $p$ are onto. A reader can easily verify that $(r, \Id): M\MonQ H \to N$ is indeed a map morphism.

When we are making a monodromy quotient of a map $M$, the new flags are orbits of a normal group 
$H \vartriangleleft \Mon(M)$. The quotienting works, because the orbits are the blocks of imprimitivity of the action $(\Flags(M), \Mon(M))$.
If we take any subgroup $K \leq \Aut(M)$ then the orbits of that subgroup are also blocks of imprimitivity for the same
action. This kind of quotients was discussed in \cite{malnicSkovieraNedela}. 
We will call such a quotient an {\em automorphism quotient}.


Having in mind the results of this section we will often say that some map is a monodromy quotient of a map $M$ if it is
isomorphic to some monodromy quotient of the map $M$.

\section{Parallel product and automorphisms}

In this section some properties of the parallel product that include the lifts of automorphisms are discussed. 

\begin{proposition}
\label{liftParallel}
If maps $M_i=(f_i, G_i, Z_i, \id_i)$, $i=1,2,$ contain  automorphisms $\alpha_W$ 
then the map $M_1\parallel M_2$ also contains the automorphism $\alpha_W$.
\end{proposition}
\begin{proof}
A parallel product is obtained as: 
$$
M := M_1 \parallel M_2  = (f_{1,2}, G := f_{1,2}(F), Z:= \Orb_{\id_{1,2}}^{G}(Z_{1,2}), \id_{1,2}).
$$ 
For a word $w \in F$, take $\alpha_w^i \in \Aut(M_i)$ and let $\alpha = (\alpha_w^1, \alpha_w^2)$. Note that in this proof
the superscripts are not the exponents but are used  as indices.
Let $z = (z_1, z_2) \in Z$ and $W = f_{1,2}(w) \in G$. Then: 
\begin{align*}
\alpha(z\cdot W) &= \alpha(z_1\cdot f_1(w), z_2\cdot f_2(w)) = (\alpha_w^{1}(z_1 \cdot f_1(w)), \alpha_w^{2}(z_2\cdot f_2(w))) =\\
&=(\alpha_w^{1}(z_1)\cdot f_1(w), \alpha_w^{2}(z_2)\cdot f_2(w)) = (\alpha_w^{1}(z_1), \alpha_w^{2}(z_2))\cdot f_{1,2}(w) =\\ 
&=\alpha(z)\cdot W.
\end{align*}
As $\alpha_w^{i}$, $i = 1, 2$, are one-to-one, $\alpha$ is indeed an automorphism. Note that $\alpha = \alpha_w \in \Aut(M)$.
\end{proof}

Since the parallel product is associative, the proposition can be generalized to a parallel product of 
a finite number of maps.

It was proven by Wilson \cite{wilsonParallel}, that if $h: M \to N$ 
is a morphism of rooted maps 
then $M \parallel N \simeq M$. This also yields $M\parallel M \simeq M$. 

\begin{proposition}
\label{uniqueParallel}
A parallel product $M \parallel N$ is the unique minimal cover over $M$ and $N$. Any cover $C$ over $M$ and $N$ is a cover
of $M \parallel N$.
\end{proposition}
\begin{proof}
Note that $(M \parallel N) \parallel C = (M \parallel C) \parallel (N \parallel C) = C \parallel C = C$.
\end{proof}

If we forget the word "rooted" in the Proposition \ref{uniqueParallel} then the proposition is not true anymore.
The example of that can be seen in Figure \ref{cyc4quo} later in Section \ref{secDecRef}.

Together with common automorphisms in two maps some other automorphisms can be present in
a parallel product. In the case where factors are re-rooted maps, the following claim 
was noted in \cite{wilsonParallel} and generalized here.

\begin{proposition}
Let $M_i = (f_i, G_i, Z_i, \id_i)$, $i = 1,\ldots, n$, be maps obtained by re-rooting a map $M$, 
$N = M_1 \parallel \cdots \parallel M_n$ be the parallel product and $\alpha$ a permutation
of components in the Cartesian product $\prod_{i=1}^{n}Z_i$. If $\alpha$ maps the orbit
of the action of the group $(f_i)_{i=1}^{n}(F)$ acting on $\prod_{i=1}^n Z_i$ containing $(\id_i)_{i = 1}^{n}$ to itself, then
$\alpha \in \Aut(M)$.  
\end{proposition}
\begin{proof}
Let $\pi \in \Syml(n)$, such that $\alpha\left((z_1, \ldots, z_n)\right) = (z_{\pi(1)}, \ldots, z_{\pi(n)})$.
Note that $f = f_1 = \ldots = f_n$, since the maps $M_i$ are obtained by re-rooting of the same map.
Let $g = (f)_{i = 1}^{n}$. For any $W \in F$, 
\begin{align*}
\alpha\left((z_1,\ldots, z_n)\cdot g(W)\right) &=
\alpha\left((z_1\cdot  f(W), \ldots, z_n\cdot  f(W))\right) = \\
&= (z_{\pi(1)}\cdot f(W), \ldots, z_{\pi(n)}\cdot f(W)) =\\
&= (z_{\pi(1)},\ldots, z_{\pi(n)})\cdot g(W) = \alpha\left((z_1, \ldots, z_n)\right)\cdot  g(W)
\end{align*}
and the result
follows.
\end{proof}

For a given map $M$ let $M^M$ denote the {\em total parallel product} of the map $M$ defined as the parallel product 
of all re-rooted maps obtained from the map $M$. 

\begin{proposition}
\label{totalInfo}
Let $M$ be an arbitrary rooted map.
\begin{enumerate}
\item
\label{prod::1}
If $M'$ and $M''$ are maps obtained from the map $M$ by re-rooting then  $\Mon(M) = \Mon(M') = \Mon(M'') = \Mon(M' \parallel M'')$.
\item
\label{prod::2}
If $M'$ and $M''$ are maps obtained from $M$ by re-rooting, such that both of them have 
a root flag in the same orbit of $\Aut(M)$, then $M'$ and $M''$ are isomorphic as rooted maps.
\item
\label{prod::3}
The total parallel product $M^M$ is a reflexible map. It is the smallest reflexible cover over the
map $M$. Any reflexible cover of $M$ is also a cover of $M^M$.
\end{enumerate}
\end{proposition}
\begin{proof}
Since $f = f_{M'} = f_{M''}$ and $f(F) \simeq (f, f)(F)$, (\ref{prod::1}) follows. 

Let $\alpha \in \Aut(M)$, 
such that $\alpha(\id_{M'}) = \id_{M''}$. Then $(\alpha, \Id)$ is an isomorphism of the rooted 
maps $M'$ and $M''$ and (\ref{prod::2}) follows.

Let $1, \ldots, n$, be the flags of the map $M$. Then $\id_{M^M} = (1, \ldots, n)$. 
Using (\ref{prod::1}) $\Mon(M^M) = \Mon(M)$. Let $W \in \Mon(M^M)$. Then 
$\id_{M^M}\cdot W = \id_{M^M}$ implies that $W$ viewed as an element of $\Mon(M)$ 
stabilizes all the flags in $M$, thus it is contained in the kernel of the action of $\Mon(M)$
on $\Flags(M)$, which is trivial. Therefore $\Mon(M^M)$ acts regularly on $\Flags(M^M)$ and thus
$M^M$ is reflexible. 

Let $N$ be any reflexible cover over $M$ and $(p,q): N \to M$ be the corresponding map morphism. 
Then $\Mon(N) = q^{-1}(\Mon(M))$. 
A reflexible map is completely determined by its monodromy group, since by Corollary \ref{isomorph} such
a reflexible map $(f, G, Z, \id)$ is isomorphic to the map $(f, G, G, 1)$, where $1 \in G$ denotes an identity element.
In our case $N = M^M$, $\Mon(N) = \Mon(M)$ and $q$ is an identity mapping.
Thus $M^M$ must be the unique minimal reflexible cover over $M$.
It is also obvious that any reflexible cover over $M$ is also a cover over $M^M$.
Therefore, (\ref{prod::3}) follows.
\end{proof}

From Proposition \ref{totalInfo}, the following corollary immediately follows.
\begin{corollary}
All re-rootings of a reflexible map are isomorphic. \qed
\end{corollary}

Therefore, when we are working with reflexible maps only, we can omit the roots, 
since any choice of root yields the same rooted map.

From the proof of Proposition \ref{totalInfo} it can be seen that the minimal reflexible cover can be obtained in a much easier 
way then by calculating $M^M$. From $M = (f, G, Z, \id)$ one just needs to construct $(f, G, G ,1)$ and this is already the 
minimal reflexible cover.

The following proposition is also very useful.
\begin{proposition}
\label{parDist}
Let $M$ and $N$ be rooted maps. 
\begin{enumerate}
\item
$\Du(M \parallel N) = \Du(M) \parallel \Du(N)$.
\item
$\Pe(M \parallel N) = \Pe(M) \parallel \Pe(N)$.
\qed
\end{enumerate}
\end{proposition}
\begin{proof}
Let $M, N$ be $(f_i, G_i, Z_i, \id_{i})$, $i = 1, 2$, respectively. Then 
\begin{align*}
\Du&(M\parallel N) = \Du\left((f_{1,2}, f_{1,2}(F), \Orb_{\id_{1,2}}^{f_{1,2}(F)}(Z_{1,2}), \id_{1,2})\right)\\
&=\left(f_{1,2}\circ\dd, (f_{1,2}\circ\dd)(F), \Orb_{\id_{1,2}}^{(f_{1,2}\circ\dd)(F)}(Z_{1,2}), \id_{1,2}\right)\\
&= \left((f_1\circ\dd, f_2\circ\dd), (f_1\circ\dd, f_2\circ\dd)(F), \Orb_{\id_{1,2}}^{(f_1\circ\dd, f_2\circ\dd)(F)}(Z_{1,2}), \id_{1,2}\right)\\
&=(f_1\circ\dd, (f_1\circ\dd)(F), Z_1, \id_1)\parallel (f_2\circ\dd, (f_2\circ\dd)(F), Z_2, \id_2)\\
&= \Du(M) \parallel \Du(N).
\end{align*}
The proof for the operation $\Pe$ is similar.
\end{proof}

For a given reflexible map $M$ we can construct a reflexible cover $N$, such that $\Du(N) = \Pe(N) = N$, 
i.e. a self-dual and a self-Petrie reflexible map. Such a map is called {\em totally symmetric}.

\begin{proposition}
Let $M$ be a reflexible map. Then the parallel product of all the maps obtained from $M$ by applying the compositions of 
the operations $\Du$ and $\Pe$ is totally symmetric and is unique minimal with these properties.
\end{proposition}
\begin{proof}
Denote by $S$ the set of all the non-isomorphic maps obtained from $M$ using the operations $\Du$ and $\Pe$ and
denote the parallel product of all the maps in $S$ by $N$. Since the parallel product operation is commutative, 
any order of the factors  in the parallel product of all the maps in $S$ always yields the map (isomorphic to) $N$. 
Let us prove that $N$ is self-dual.
Since in the set $S$ there are all non-isomorphic maps obtained by the operations $\Du$ and $\Pe$ 
and the operations are involutions, 
performing $\Du$ on all the elements of $S$ yields the same set. Therefore the parallel product yields a map isomorphic
to $N$ and by Proposition \ref{parDist} the map $N$ is self-dual.
Similarly we show that $N$ is self-Petrie. 
If $N'$ is a cover of $M$ it follows that $\Du(N')$ must be a cover of $\Du(M)$ (and similarly for the operation $\Pe$).
Using that, a reader can easily verify the minimality and the uniqueness.
\end{proof}

Some of those properties of the parallel product were already noted in
\cite{wilsonParallel} without a proof. 
This theory can be extended in several directions. A possible extensions include abstract polytopes \cite{schulteMullen}. 
Using the constuctions in this section one can extend the results to abstract polytopes and get similar 
results to the ones by Hartley \cite{hartley}.

\section{A parallel-product decomposition of a map}

We consider factorizing a map $M$ as a parallel product. 
The factors are always 
the images of map morphisms. Our aim is to find criteria for splitting the map as a parallel product 
of two maps which are monodromy quotients. Monodromy quotients are of special interest, 
because all the automorphisms $\alpha_W \in M$ project. In particular, a monodromy quotient
of a reflexible map is reflexible.

A map $M$ is {\em parallel-product decomposable} if it is a non-trivial parallel product of two maps, such that 
the two maps are monodromy quotients of $M$.

\begin{theorem} (Decomposition theorem) A map $M = (f, G, Z, \id)$  is parallel-product decomposable
if and only if there are at least two normal non-transitive subgroups $\Ha, \Hb \vartriangleleft G$, such that 
$\Ha \cap \Hb = \{1\}$ and $G_{\id}\Ha \cap G_{\id}\Hb = G_{\id}$.
\label{decomposition}
\end{theorem}
\begin{proof}
Let $M = M_1 \parallel M_2 = (f_{1,2},~ G:=f_{1,2}(F),~ Z:=\Orb_{\id_{1,2}}^{G}(Z_{1,2}),~ \id := \id_{1,2})$  
be a non-trivial parallel product of maps,
where $M_i = (f_i, \Gi, $ $Z_i, \id_i)$, $i = 1, 2$. Note that the indices in the names of the groups are written as superscripts since
subscripts are used for denoting stabilizers. The coordinate projections
$(p_i, q_i): Z\times G \to Z_i \times \Gi$ are the covering projections of the maps $M \to M_i$.
Denote the kernels of the epimorphisms $q_i$ by $\Hi$. These are normal subgroups in $G$ 
and $\Ha \cap \Hb = \{(1, 1)\}$. Since the factors of the parallel product are monodromy quotients, they must 
be the monodromy quotients by these two normal subgroups.
For monodromy quotients it is true:  $q_i^{-1}(\Gi_{\id_i}) = G_{\id}\Hi$.
But since $G_{\id} = q_1^{-1}(\Ga_{\id_1}) \cap q_2^{-1}(\Gb_{\id_2})$ it immediately follows:  
$G_{\id}\Ha \cap G_{\id}\Hb = G_{\id}$.
Thus if $M$ is a nontrivial parallel product of two maps that are the monodromy quotients of the product 
then it meets the conditions of the theorem.

Now, let $M = (f, G, Z, \id)$. By Corollary \ref{isomorph} we can assume that
$M = (f, G,$ $G/G_{\id}, G_{\id})$. 
Let $\Ha$, $\Hb$ be the normal subgroups meeting the conditions of the theorem.
A trivial parallel-product decomposition would be obtained if one of the factors would be isomorphic to $M$ or to 
the trivial map. In the first case this would mean $\Hi \leq G_{\id}$, but since the action $(G/G_{\id},G)$ is faithful
this cannot happen. The second case is prevented by the non-transitivity condition.

By Proposition \ref{monQuoDef} the monodromy quotients of $M$ by $\Hi$, $i = 1, 2$, are isomorphic to the maps 
$M_i := (f_i, \Gi, Z_i, \id_i)$, where $\Gi := G/\Hi$, 
$f_i := q_i\circ f$, $q_i: G \to G/\Hi$ is a natural epimorphism, $Z_i := G/G_{\id}\Hi$ and $\id_i = G_{\id}\Hi$. Let $(p_i, q_i): M \to M_i$ be the corresponding
covering projections as in Proposition \ref{monQuoDef}. It is easy to see that 
$p_i: G/G_{\id} \to G/G_{\id}\Hi$ is defined by $p: G_{\id}g \to G_{\id}\Hi g$, for any $g \in G$.

Let   $M_1 \parallel M_2 = (f_{1,2}, K, X, \idd)$, where 
$K = f_{1,2}(F)$, $\idd = (G_{\id}\Ha,G_{\id}\Hb)$ 
and $X$ is an orbit of the naturally induced action $(G_{\id}\Ha \times  G_{\id}\Hb, K)$ containing
$\idd$.
We will show that $M_1 \parallel M_2$ is isomorphic to $M$ and therefore we have to find an
isomorphism $\psi: K \to G$ and a bijection $\phi: X \to G/G_{\id}$, such that $(\phi, \psi): M_1 \parallel M_2 \to M$ is 
a map isomorphism.

Let $W \in K$. Then there exists $w_1 \in F$, such that $f_{1,2}(w_1) = W$. Let $\psi(W) := f(w_1)$.  
First we verify that $\psi$ is well defined. It is true that  $f_{1,2}(w_1) = (q_1\circ f(w_1), q_2 \circ f(w_1))$.
If there is some other $w_2 \in F$, such that $f_{1,2}(w_2) = W$, we get 
$q_i\circ f(w_1) = q_i\circ f(w_2)$, $i \in \{1, 2\}$.  This means $f(w_1)f^{-1}(w_2) \in \Ha \cap \Hb = \{1\}$ 
and it follows $f(w_1) = f(w_2)$. Hence, the mapping $\psi$ is well defined. Now we have to see that 
$\psi$ is a homomorphism of the groups. Let $g = (g_1, g_2), h = (h_1, h_2) \in K$. There are $w_1, w_2 \in F$,
such that $g = f_{1,2}(w_1)$ and $h = f_{1,2}(w_2)$. Then $\psi(g) = f(w_1)$ and 
$\psi(h) = f(w_2)$.
Since
\begin{align*}
f_{1,2}(w_1 w_2) &= (f_1(w_1 w_2), f_2(w_1 w_2)) = (f_1(w_1)f_1(w_2), f_2(w_1)f_2(w_2))\\
&= (g_1 h_1, g_2 h_2) = gh,
\end{align*}
then $\psi(gh) = f(w_1 w_2) = f(w_1) f(w_2) = \psi(g)\psi(h)$ and $\psi$ is a homomorphism. Obviously, 
it is an epimorphism. Let $g \in \ker \psi$. 
There exists $w \in F$, such that $f_{1,2}(w) = g$ and $\psi(g) = f(w) = 1$.
Thus $f_1(w) = f_2(w) = 1$ and since $g = (f_1(w), f_2(w))$, it follows that $g = 1$ and $\psi$ must be an isomorphism.

Let $z \in X$. Then 
$z = (G_{\id}\Ha f_1(w), G_{\id}\Hb f_2(w))$, for some $w \in F$. 
Define $\phi: X \to G/G_{\id}$ by $\phi: z \mapsto G_{\id}f(w)$.
There may exist another $w' \in F$, such that $z = (G_{\id}\Ha f_1(w'), G_{\id}\Hb f_2(w'))$.
Then 
\begin{align*}
(G_{\id}\Ha f_1(w')f_1^{-1}(w), G_{\id}\Hb f_2(w')f_2^{-1}(w)) &= (G_{\id}\Ha,G_{\id}\Hb)f(w')f^{-1}(w)\\ 
                            &= (G_{\id}\Ha,G_{\id}\Hb).
\end{align*}
Thus by the assumption of the theorem $f(w')f^{-1}(w) \in G_{\id}\Ha \cap G_{\id}\Hb = G_{\id}$
and $\phi$ is well defined. Similarly we can see that $\phi$ is one-to-one.
Since for any $w \in F$, it follows $(G_{\id}\Ha f_1(w), G_{\id}\Hb f_2(w)) \in X$, the mapping $\phi$
is onto.

Now we will verify that $(\phi, \psi)$ is an isomorphism of the maps $M_1 \parallel M_2$ and $M$. 
Obviously, $\phi(\idd) = G_{\id}$ and  $\psi \circ f_{1,2} = f$.
Let $g \in K$. Then there exists $w_1 \in F$, such that $g = f_{1,2}(w_1)$ and $\psi(g) = f(w_1)$.
Let $z = (G_{\id}\Ha f_1(w_2), G_{\id}\Hb f_2(w_2)) \in X$ for some $w_2 \in F$. Then
\begin{align*}
zg &= (G_{\id}\Ha f_1(w_2), G_{\id}\Hb f_2(w_2))\cdot f_{1,2}(w_1)\\
&= (G_{\id}\Ha f_1(w_2 w_1), G_{\id}\Hb f_2(w_2 w_1))
\end{align*}
and $\phi(zg) = G_{\id}f(w_2 w_1) = G_{\id}f(w_2)f(w_1)$. 
Also, $\phi(z)\psi(g) = G_{\id}f(w_2)f(w_1)$. Therefore, a pair $(\phi, \psi)$ is an isomorphism.
\end{proof}

\section{Degeneracy of reflexible maps}

In this section reflexible maps are classified into three families according to their degeneracy. 
The classification will be used in the following section, where all parallel-product indecomposable degenerate maps will be presented.

Let $M$ be a reflexible map with a presentation of the monodromy group of the form
$$
\Mon(M) = \langle \TT, \LL, \RR ~\vert~ W_1^{e_1}= W_2^{e_2} = \ldots = W_k^{e_k}  = 1\rangle,\quad  e_i \geq 1, k \geq 7,   
$$
such that $W_1 = \TT, W_2 = \LL, W_3 = \RR, W_4 = \LL\TT$, $W_5 = \RR\TT$, $W_6 = \RR\LL$, $W_7 = \TT\LL\RR$,
where  
$e_1, e_2, e_3, e_4 \in \{1,2\}$, and where $W_i, i \geq 8$ are words in $\Mon(M)$,
such that the group is finite. 
Also, all $e_i$ are the actual orders of the corresponding elements (words). 
The set of words $\{W_1, \ldots, W_k\}$ is called a {\em context}. Any context 
contains at least the words $W_1, \ldots, W_7$.
In the context chosen, a
monodromy group can be denoted by a vector 
$\Mon(M) = (e_1, e_2, \ldots, e_k)$ or $\Mon(M) = (e_i)_{i = 1}^k$.
When for a given map $M$ the words in the context $C$ are sufficient to define $\Mon(M)$, the context is said to be 
{\em sufficient}.
A monodromy group $\Mon(M)$ can be easily obtained from the vector and 
the obtained reflexible map is $M = (f, \Mon(M), \Mon(M), 1)$, where $f$ is a homomorphism mapping $t \mapsto W_1$,
$l \mapsto W_2$, $r \mapsto W_3$ and $1 \in \Mon(M)$. Sometimes the notation
is abused and the map is denoted directly by the corresponding vector.
It is obvious that any monodromy group of a reflexible map can be written in the form described above,
but some of the maps need larger contexts (i.e. more words $W_i$, $i \geq 8$).

For two contexts $C_1$ and $C_2$ the {\em common context} is $C_1 \cup C_2$. Obviously, if some map
is represented in a context $C_1$, it can be also easily represented  in $C_1 \cup C_2$ by calculating
the orders of the words in $C_2\setminus C_1$ and adding those (redundant) relations.

A map $M$
is {\em slightly-degenerate} if in any sufficient context $C$ it follows $e_i \geq 2$, for all $i = 1,\ldots,7$, and at least one 
of $e_5, e_6, e_7$ equals to 2.
It is {\em degenerate} if at least one of $e_i$, $i = 1, \ldots, 7$, equals to 1. If a map is not degenerate or 
slightly-degenerate then it is {\em non-degenerate}. In this case $e_i \geq 3$, $i = 5,6, 7$.

Note that in any sufficient context of a map $M$ the words $W_i$, $i = 1, \ldots, 7$, are exactly the generators and the relations that determine
the map's properties, such as the degrees of the vertices, the co-degrees of the faces and the sizes of the Petrie circuits.

\begin{lemma}
\label{vectorProduct}
Let $M = (e_i)_{i = 1}^k$, $N = (f_i)_{i=1}^k$ be two reflexible maps represented in the common context. Then 
$M\parallel N = (\lcm(e_i,f_i))_{i = 1}^k$.
\end{lemma}
\begin{proof}
We can view both groups $\Mon(M)$ and $\Mon(N)$ as quotients of a free group $F_0 := \langle \TT, \LL, \RR \rangle$. 
Let $\Ha$ be the normal closure in $F_0$ of the set $\{ W_i^{e_i}\}_{i=1}^{k}$ 
and $\Hb$ be the normal closure in $F_0$ of $\{W_i^{f_i}\}_{i = 1}^{k}$. Then $\Mon(M) = F_0/\Ha$ and $\Mon(N) = F_0/\Hb$.
Let $g$ be an element of the intersection $\Ha \cap \Hb$. Then $g$ can be expressed as a finite product of 
conjugates and powers of conjugates of 
$W_i$. Since everything is happening in the free group, any exponent of $W_i$ in the expression of $g$ must
be divisible by $e_i$ and $f_i$ and thus by $\lcm(e_i, f_i)$. Thus $\Ha \cap \Hb$ is exactly the normal closure of
of $\{W_i^{\lcm(e_i,f_i)}\}_{i = 1}^{k}$ and this set determines the relations of $F_0/(\Ha \cap \Hb)$ in the finite presentation. 
Let $f_i: F_0 \to F_0/\Hi$, $i = 1, 2$,
be the natural quotient projections. Let $f := (f_1, f_2)$. Then $\ker f = \Ha \cap \Hb$ and $F_0/(\Ha\cap \Hb) \simeq f(F_0)$.
But $M = (f_1, F_0/\Ha, F_0/\Ha, \Ha)$ and $N = (f_2, F_0/\Hb, F_0/\Hb, \Hb)$ 
and thus $F_0/(\Ha\cap \Hb) = \Mon(M \parallel N)$.

\end{proof}

When  analyzing the existence of some class of reflexible maps for which the context $C = \{W_1, \ldots, W_7\}$
is sufficient,
we can use triality. Note
that the operations $\Du$ and $\Pe$ permute the triple $(e_1, e_2, e_4)$ with the same permutation 
as the triple $(e_5, e_6, e_7)$. To describe the action of $\Du$ and $\Pe$ on the indices $i = 1, \ldots, 7$ 
of $e_i$, we can represent $\Du$ as a permutation $(1,2)(5,6)$ and $\Pe$ as $(2,4)(6,7)$.

\begin{proposition}
All degenerate reflexible maps are shown in Table \ref{degenerate}.
\end{proposition}
\noindent
\begin{table}[ht]
\centering
\bigskip
\begin{tabular}{llccccccccc}

{\bf Name} 	& ($\TT$, &$\LL$, &$\RR$, &$\TT\LL$, &$\TT\RR$, &$\LL\RR$, &$\TT\LL\RR$) & $|Mon(M)|$\\
$\DM_1$         & ($1$,&$1$,&$1$,&$1$,&$1$,&$1$,&$1$)	& 1		\\  
$\DM_2$         & ($1$,&$1$,&$2$,&$1$,&$2$,&$2$,&$2$)	& 2		\\
$\DM_3$         & ($2$,&$1$,&$1$,&$2$,&$2$,&$1$,&$2$)	& 2		\\
$\DM_4$         & ($1$,&$2$,&$1$,&$2$,&$1$,&$2$,&$2$)	& 2		\\
$\DM_5$         & ($2$,&$2$,&$1$,&$1$,&$2$,&$2$,&$1$)	& 2		\\
$\DM_6(k)$, $k > 0$      & ($2$,&$1$,&$2$,&$2$,&$k$,&$2$,&$k$)	& $2k$		\\
$\DM_7(k)$, $k > 0$      & ($1$,&$2$,&$2$,&$2$,&$2$,&$k$,&$k$)	& $2k$		\\
$\DM_8(k)$, $k > 0$	& ($2$,&$2$,&$2$,&$1$,&$k$,&$k$,&$2$)     & $2k$	\\
$\DM_9$         & ($2$,&$2$,&$1$,&$2$,&$2$,&$2$,&$2$)	& 4	\\
$\DM_{10}$	& ($2$,&$2$,&$2$,&$2$,&$1$,&$2$,&$2$)	& $4$	\\
$\DM_{11}$	& ($2$,&$2$,&$2$,&$2$,&$2$,&$1$,&$2$)	& $4$	\\
$\DM_{12}$	& ($2$,&$2$,&$2$,&$2$,&$2$,&$2$,&$1$)	& $4$
\end{tabular}
\caption{\label{degenerate}Degenerate reflexible maps.
}

\end{table}

\begin{proof}
First we  prove that all the monodromy groups in Table \ref{degenerate} are uniquely determined by the context
$C = \{W_1, \ldots, W_7\}$. For all the maps in the table except $\DM_{i}(k)$, $i = 6,7,8$, this is pretty obvious. 
By triality it is enough to check the group of $\DM_6(k)$. The relations here determine 
a dihedral group $\Dih_{2k}$ generated
by $\TT$ and $\TT\RR$ that commute. 
One can easily see that any quotient of $\Dih_{2k}$ strictly decreases the orders of at least one of the 
(projected) generators.

Now we will make an analysis of what kind of degenerate maps can occur.
Let $e_1 = e_2 = 1$. Then $e_4 = 1$. If $e_3 = 1$  
we get $\DM_1$. If $e_3 = 2$ then it must be 
$e_5 = e_6 = e_7 = 2$ ($\DM_2$).
Now, let $e_1 = 1$ and $e_2 = 2$. Since $e_4 = 1$ implies $e_2 = e_1$, it must be $e_4 = 2$.
If $e_3 = 1$ then it must be $e_5 = 1$, $e_6 = e_7 = 2$ ($\DM_4$ and by triality $\DM_3$ and $\DM_5$).
If $e_3 = 2$ then $e_5 = 2$ and $e_6 = e_7 = k \geq 1$ ($\DM_7(k)$ and by triality $\DM_6(k)$ and $\DM_8(k)$).
By triality, all the possibilities where one of $e_1, e_2, e_4$ is 1 are exhausted.
Assume  $e_1 = e_2 = e_4 = 2$. If $e_3 = 1$ then $e_5 = e_6 = e_7 = 2$ ($\DM_9$).
Let now $e_3 = 2$. Since map has to be degenerate, one of $e_5, e_6, e_7$ must be equal to 1.
By triality we can assume $e_5 = 1$. Then it must be $e_6 = e_7 = 2$, otherwise the orders 
$e_1, e_2$ collapse ($\DM_{10}$, $\DM_{11}$, $\DM_{12}$). 
This exhausts all the possibilities for degenerate maps.
\end{proof}	       
A similar analysis of degenerate maps was done in \cite{caihengli}, but \v Sir\'a\v n's definition of degeneracy is
different from ours. By \v Sir\'a\v n, a reflexible map $M$ is degenerate 
if one of the generators $x = \alpha_{\LL}$, $y = \alpha_{\TT}$, $z= \alpha_{\RR} \in \Aut(M)$ equals to the identity.   
It is easy to see that
\v Sir\'a\v n's degeneracy is equivalent to saying that one of $e_1$, $e_2$ or $e_3$ is equal to 1.
Unfortunately, in \cite{caihengli} they  forgot to include the map $\DM_5$.
They also use similar names for degenerate maps. Thus their maps $\DM_1$, $\ldots$, $\DM_7$ correspond 
to ours $\DM_1$, $\DM_2$, $\DM_4$, $\DM_3$, $\DM_6$, $\DM_7$ and 
$\DM_9$, respectively.

\begin{figure}[htp]
\centering
\centerline{\epsfig{file=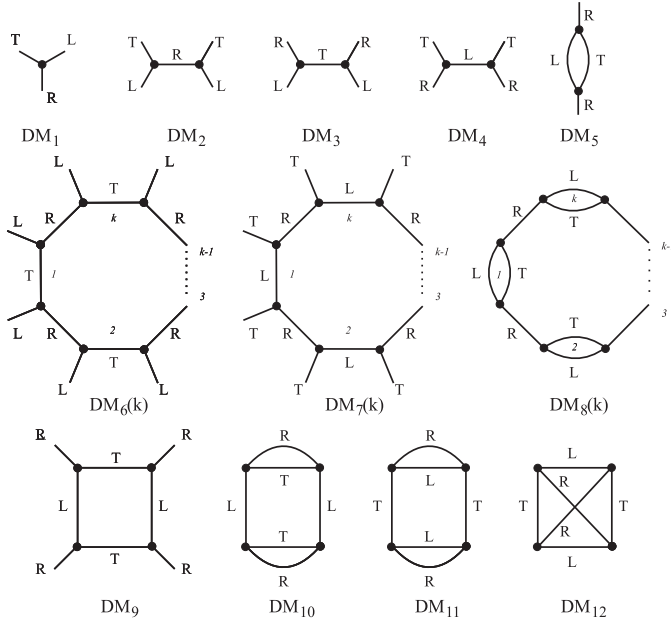,width=0.7\textwidth}}
\caption{\label{degenerateCo} Flag graphs of degenerate reflexible maps.}
\end{figure}

In Figure \ref{degenerateCo} all the flag graphs for degenerate maps are shown.

If a reflexible map is not degenerate then all the involutions 
$\TT$, $\LL$, $\RR$, $\TT\LL$ are fixed-point-free.
Such a map corresponds to a reflexible 2-cell embedding of some graph into a
compact closed surface.
Slightly-degenerate maps can be constructed using the operations $\Du$ and $\Pe$ 
from a reflexible embedding of a cycle in some compact closed surface. 
The only possible such 2-cell embeddings are the embeddings of $k$-cycle in the sphere, denoted 
by $\varepsilon_k$, and 
in the projective plane with the $k$-cycle embedded as a non-contractible curve, denoted by $\delta_k$.
Here the names are adopted from \cite{wilsonBic}.

The monodromy group presentations of maps $\varepsilon_k$ and $\delta_k$ are shown in Table \ref{slightlyDegenerate}.

\noindent
\begin{table}[ht]
\centering
\bigskip
\begin{tabular}{llc}
{\bf Name} 				& {\bf Additional relations} & {\bf Order}\\
$\varepsilon_k$, $k > 0$ even 		& $(\LL\RR)^k, (\TT\LL\RR)^k$     	&$4k$\\
$\varepsilon_k$, $k > 1$ odd 		& $(\LL\RR)^k, (\TT\LL\RR)^{2k}$ 	&$4k$\\
$\delta_k$, $k > 0$ even 		& $\TT(\LL\RR)^k, \TT(\TT\LL\RR)^{k}$   &$4k$\\
$\delta_k$, $k > 1$ odd 		& $(\LL\RR)^{2k}, (\TT\LL\RR)^{k}$      &$4k$\\
\end{tabular}
\caption{\label{slightlyDegenerate}
A monodromy group of each map in this table is obtained as 
$\langle \TT, \LL, \RR~|~ \TT^2 = \LL^2 = \RR^2 = (\TT\LL)^2 = (\RR\TT)^2 = \ldots = 1\rangle$, where instead
of ''\ldots'' one should put the additional relations. 
All slightly-degenerate reflexible maps can be constructed from the
maps in this table by using the operations $\Du$ and $\Pe$.
Note that $\varepsilon_1 = \DM_{11}$ and 
$\delta_1 = \DM_{12}$ and thus degenerate and not included in Table \ref{slightlyDegenerate}.
}
\end{table}

\section{Parallel-product decomposition of reflexible maps}
\label{secDecRef}
For reflexible maps the decomposition theorem (Theorem \ref{decomposition}) can be more specialized.

\begin{theorem}
\label{reflexibleDecomposition}
A reflexible map $M$ is parallel-product decomposable if and only if $\Mon(M)$ (and therefore also $\Aut(M)$) 
contains at least two non-trivial 
minimal normal subgroups. 
\end{theorem}
\begin{proof}
Since the monodromy group of a reflexible map is regular, the stabilizer is trivial.  The conditions of 
Theorem \ref{decomposition} are reduced to the existence of two non-trivial normal subgroups $\Ha$ and $\Hb$, such
that $\Ha \cap \Hb = \{1\}$. But in a finite group such subgroups exist if and only if two minimal non-trivial 
normal subgroups exist. Since for reflexible maps the monodromy group is isomorphic to the automorphism group, the result
follows.
\end{proof}

\begin{example}
To demonstrate how quotienting and the parallel-product decomposition work, see the examples in Figures \ref{cyc4mon} and 
\ref{cyc4quo}. In both figures the map $M$ we are quotienting is a 4-cycle on the sphere.
In Figure \ref{cyc4mon}, $M$ and its quotients are represented by flag graphs. 
We note that the monodromy group $\Mon(M)$ is isomorphic to the group $\ZZ_2 \times \Dih_4$. This 
group has exactly 3 minimal normal subgroups. The flag graphs of each of the corresponding monodromy quotients are shown.
All these maps are reflexible. By Theorem \ref{reflexibleDecomposition}
a parallel product of any two yields the original map $M$.

\begin{figure}[htp]
\centering
\centerline{\epsfig{file=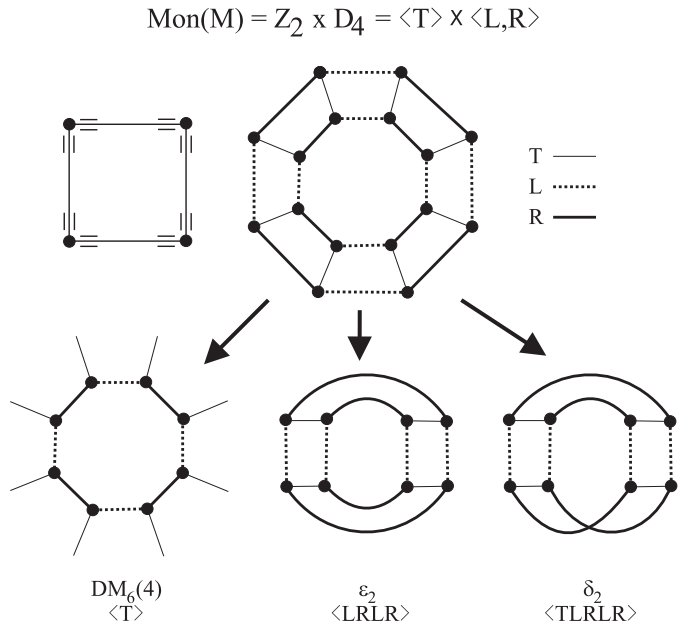,width=0.7\textwidth}}
\caption{\label{cyc4mon} Monodromy quotients of $C_4$ on the sphere that yield a non-trivial parallel-product decomposition.}
\end{figure}

In Figure \ref{cyc4quo}, a different quotient is obtained. The quotient arises as an automorphism quotient 
from the orbits of the automorphism that
rotates the flags around the vertex in the lower left corner. 
The obtained map is not reflexible. One can easily see that in the quotient there are 2 orbits of the automorphism group on the flags, 
namely the orbit of the flags around the vertices of degree 1 and the orbit of the flags around the vertex of degree 2. 
If we re-root the maps in a way, such that the root flags are in the different orbits and make a parallel product of them, 
we obtain the smallest 
reflexible cover (by Proposition \ref{totalInfo}) which is again the map $M$. 
Note that the monodromy groups of $M$ and its quotient are isomorphic.

\begin{figure}[htp]
\centering
\centerline{\epsfig{file=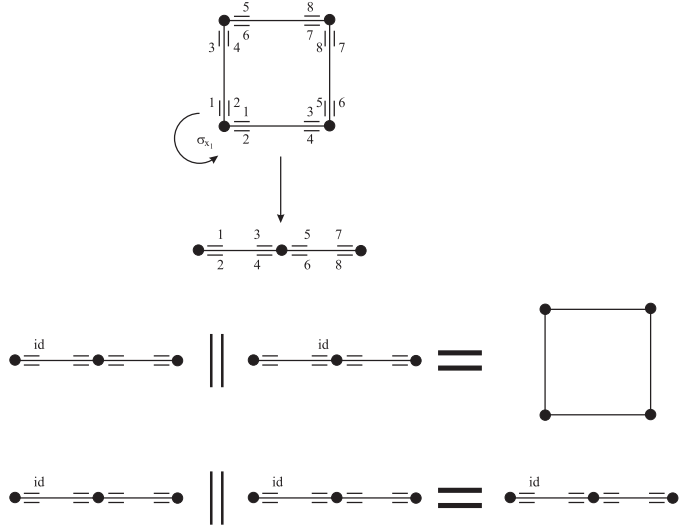,width=0.7\textwidth}}
\caption{\label{cyc4quo} An automorphism quotient of $C_4$ and a demonstration of Proposition \ref{totalInfo}.}
\end{figure}
\end{example}

Thus a parallel-product indecomposable reflexible map is any reflexible map $M$, 
such that either $\Mon(M)$ is  a simple group or 
$\Mon(M)$ has a unique minimal normal subgroup. The latter groups are called {\em monolithic groups} and 
the unique minimal normal subgroup is called a {\em monolith}.
Since the operations $\Du$ and $\Pe$ preserve a monodromy group, the operations are invariant for
the parallel-product indecomposability.

\begin{proposition}
The map $\DM_6(k)$ ($\DM_7(k)$, $\DM_8(k)$), $k > 2$ is parallel-product decomposable if and only if
$k$ is not a prime power. 
\end{proposition}
\begin{proof}
Number $k$ is not a prime power if and only if there exist $a,b > 1$, such that $\gcd(a,b) = 1$ and 
$k = ab$.
Using Lemma \ref{vectorProduct} and Table \ref{degenerate} it is easy to see that for any 
$a,b > 1$, 
$\DM_6(a)\parallel \DM_6(b) \simeq \DM_6(\lcm(a, b))$.
Nontrivial factors of $\DM_6(k)$ can be only degenerate maps with
$\LL = 1$, so only: $\DM_6(l)$, $l \geq 1$, $\DM_2$ and $\DM_3$. Since $\DM_2$  and $\DM_3$ are quotients of 
any $\DM_6(l)$, $l > 2$, a parallel product with $\DM_6(l)$ absorbs them. 
Also $\DM_2 \parallel \DM_3 \simeq $ $\DM_2 \parallel \DM_6(1) \simeq$ $\DM_3 \parallel \DM_6(1) \simeq$ $\DM_6(2)$.
So if $k > 2$ and $\DM_6(k)$ is parallel-product decomposable, then it must be a product of two 
factors of the form
$\DM_6(l)$. By Table \ref{degenerate} and Proposition \ref{vectorProduct} this is possible only when 
the conditions of the lemma are fulfilled.
Using triality, the proofs for $\DM_7(k)$ and $\DM_8(k)$ immediately follow.
\end{proof}

The monodromy groups of the maps $\DM_i$, $i = 9,10,11,12$, are isomorphic to $\ZZ_2 \times \ZZ_2$ and thus by 
Theorem \ref{reflexibleDecomposition} the maps are parallel-product decomposable. The monodromy groups of $\DM_i$, $i = 1,2,3,4,5$, 
are either trivial or isomorphic to
$\ZZ_2$, implying that those maps are parallel-product indecomposable.

The following corollary immediately follows. 

\begin{corollary}
All degenerate reflexible maps are parallel-product indecomposable except:
\begin{enumerate}
\item
$\DM_5(k)$, $\DM_6(k)$ and $\DM_7(k)$, for $k = 2$ and any $k > 2$ which is not a power of a prime,
\item
$\DM_{9}$, $\DM_{10}$, $\DM_{11}$ and
$\DM_{12}$.
\qed
\end{enumerate}
\end{corollary}

\begin{proposition}
The only parallel-product indecomposable slightly-degenerate maps are the maps $\delta_k$, where $k = 2^n$, $n \geq 1$.
\end{proposition}
\begin{proof}
Since $\Pe(\varepsilon_k) \simeq \delta_k$, for $k$ odd, we have to consider only the parallel-product decompostions of 
maps $\varepsilon_k$  for all $k > 1$ and $\delta_k$, for $k > 1$ even.

Take a context $C = \{\TT, \LL, \RR, \TT\LL, \RR\TT, \RR\TT, \TT\LL\RR \}$. In this context
$\varepsilon_k = (2,2,2,$ $2,2,k,k)$, for $k > 0$ even,  $\varepsilon_k = (2,2,2,2,2,k,2k)$, for $k > 1$ odd,  
$\DM_{3} = (2,1,1,2,2,1,2)$ and $\DM_7(k) = (1,2,2,2,2,k,k)$. By Proposition \ref{vectorProduct}, it follows 
$\varepsilon_k \simeq \DM_7(k)\parallel \DM_3$, for any $k > 1$. 

Now, let $k > 0$ and let $l \geq 1$ be any odd number. We will prove that 
$\delta_{2^k l} \simeq \DM_7(2^k l) \parallel \delta_{2^k}$.
This would mean that for any even $u$ not equal to the power of 2, $\delta_{u}$ is parallel-product decomposable.
The monodromy groups for $\delta_{2^k}$ and $\DM_7(2^{k}l)$ are defined by relations: 
\begin{align*}
&\Mon(\delta_{2^k}): \TT^2=\LL^2=\RR^2=(\TT\LL)^2=(\RR\TT)^2 = 1, (\RR\LL)^{2^k} = (\TT\LL\RR)^{2^{k}} = T,\\
&\Mon(\DM_7(2^{k}l)): \TT=\LL^2=\RR^2=(\TT\LL)^2=(\RR\TT)^2 = (\RR\LL)^{2^k l} = (\TT\LL\RR)^{2^{k} l} = 1.
\end{align*}
Hence a monodromy group of the parallel product is defined by relations

\begin{align*}
&\TT^2=\LL^2=\RR^2=(\TT\LL)^2=(\RR\TT)^2 = 1, (\RR\LL)^{2^k l} = (\TT\LL\RR)^{2^{k}l} = T
\end{align*}
and thus congruent to the monodromy group of a map $\delta_{2^k l}$.

For a given map $M$, denote by $e_5(M)$, $e_6(M)$ and $e_7(M)$ the exponents of the words $\RR\TT$, $\RR\LL$, $\TT\LL\RR$,
respectively. For $\delta_{2^n}$ it follows $e_5 = 2$, $e_6 = e_7 = 2^{n+1}$. Since these values are powers of 2 and $\lcm(2^x, 2^{y}) = \max(2^x, 2^{y})$, 
at least one of $e_5$, $e_6$, $e_7$ must be reached with the corresponding  values $e_5', e_6', e_7'$  and
$e_5'', e_6'', e_7''$ in two possible factors. 
Therefore, one of the factors should be one of $\DM_7(2^{n+1})$, $\delta_{2^{n}}$ or $\varepsilon_{2^{n+1}}$. 
In the case of $\delta_{2^{n}}$, we would not get a non-trivial product. In the case of 
$\varepsilon_{2^{n+1}}$ the parallel product would be orientable, while $\delta_{2^{n}}$ is not orientable.
Therefore, if we have a parallel-product decomposition, one of the factors must be $\DM_7(2^{n+1})$. 
Then the other factor cannot be a degenerate map, because the context $C$ is not sufficient to obtain
the map $\delta_{2^n}$. Hence, one of the factors must be a map $\delta_l$, for some $l = 2^u$, $u < n$.
But one can easily verify that in this case $\DM_7(2^{n+1}) \parallel \delta_l \simeq \varepsilon_{2^{n+1}}$.
Thus $\delta_{2^n}$, $n \geq 1$ is parallel-product indecomposable.

\end{proof}

Using computer programs {\sc Lowx} \cite{smallgenus} and {\sc Magma} \cite{magma} all non-degenerate 
reflexible maps were calculated up to 100 edges. The results of the calculation match with 
{\em Wilson's census of rotary maps} \cite{census}. Among them, the ones with the monolithic monodromy group
were selecteed and they are shown in Table \ref{minimalNonDeg}.

\begin{theorem}
Up to triality, all parallel-product indecomposable non-degenerate reflexible maps up to 100 edges are presented in Table \ref{minimalNonDeg}.
\qed
\end{theorem}

\begin{center}
\begin{longtable}{|l|c|ccc|l|l|}
\caption[]{Parallel-product indecomposable non-degenerate reflexible maps up to triality and up to 100 edges.
A presentation of any of the corresponding monodromy groups can be obtained by using a presentation
$\langle\TT, \LL, \RR~|~\TT^2 = \LL^2 = \RR^2 = (\TT\LL)^2 = (\RR\TT)^{e_5} = (\RR\LL)^{e_6} = (\TT\LL\RR)^{e_7} = \ldots = 1\rangle$,
where the corresponding additional relations should be put instead of ''\ldots''.
}\label{minimalNonDeg} \\

\hline \multicolumn{1}{|l|}{\textbf{Name}} & 
\multicolumn{1}{c|}{\textbf{$|\Mon|$}} & 
\multicolumn{1}{c|}{\textbf{$e_5$}} &
\multicolumn{1}{c|}{\textbf{$e_6$}} &
\multicolumn{1}{c|}{\textbf{$e_7$}} &
\multicolumn{1}{c|}{\textbf{Additional relations}} &
\multicolumn{1}{c|}{\textbf{Monolith}} 
\\ \hline 
\endfirsthead

\hline \multicolumn{1}{|l|}{\textbf{Name}} & 
\multicolumn{1}{c|}{\textbf{$|\Mon|$}} & 
\multicolumn{1}{c|}{\textbf{$e_5$}} &
\multicolumn{1}{c|}{\textbf{$e_6$}} &
\multicolumn{1}{c|}{\textbf{$e_7$}} &
\multicolumn{1}{c|}{\textbf{Additional relations}} &
\multicolumn{1}{c|}{\textbf{Monolith}} 
\\ \hline 
\endhead

\hline \multicolumn{7}{|r|}{{Continued on next page}} \\ \hline
\endfoot

\hline \hline
\endlastfoot

$\MN_{1}$  &24 &3 &3  &4   &  &															$\ZZ_2^2$\\
$\MN_{2}$  &32 &4 &8  &8   & $(\RR\TT\RR\LL)^2$,  $(\LL\RR\TT)^2(\LL\RR)^2$ 									&$\ZZ_2$\\
$\MN_{3}$  &60 &3 &5  &5   &															&$\{1\} \leq A_5$\\
$\MN_{4}$  &64 &4 &4  &4   &  								 							&$\ZZ_2$\\
$\MN_{5}$  &64 &4 &8  &8   &   $((\LL\RR)^2\TT)^2$        											&$\ZZ_2$\\
$\MN_{6}$  &64 &4 &16 &16  &  $(\RR\TT\RR\LL)^2$, $\TT\LL\RR\TT(\LL\RR)^7$   									&$\ZZ_2$\\
$\MN_{7}$  &72 &4 &4  &6   &  $(\RR\TT\RR\LL)^3$												&$\ZZ_3^2$\\
$\MN_{8}$  &96 &3 &8  &12  &  $(\LL\RR\TT\LL\RR)^2\TT(\LL\RR)^2\TT$       									&$\ZZ_2$\\
$\MN_{9}$  &96 &6 &8  &12  & $\TT\LL(\RR\TT)^2(\LL\RR)^3$    											&$\ZZ_2$\\
$\MN_{10}$ &108&3 &6  &6   & 															&$\ZZ_3$\\
$\MN_{11}$ &120&4 &5  &6   & $\LL(\RR\TT\RR\LL)^2(\RR\TT)^2$											&$A_5 \leq S_5$\\
$\MN_{12}$ &120&6 &6  &6   & $\LL(\RR\TT)^2\RR\LL(\RR\TT)^3$, $\TT(\LL\RR)^3\TT\RR(\LL\RR)^2$							&$A_5 \leq S_5$\\                                                    
$\MN_{13}$ &128&4 &4  &8   & $(\RR\TT\RR\LL)^4$													&$\ZZ_2$\\
$\MN_{14}$ &128&4 &16 &16  & $(\LL\RR\TT)^2(\RR\LL)^2(\RR\TT)^2$, $(\LL\RR)^2\TT(\LL\RR)^6\TT$  						&$\ZZ_2$\\
$\MN_{15}$ &128&4 &32 &32  &  $(\RR\TT\RR\LL)^2$, $(\LL\RR\TT)^2(\LL\RR)^{14}$									&$\ZZ_2$\\
$\MN_{16}$ &128&8 &8  &8   & $(\LL\RR\TT)^2(\RR\LL)^2(\RR\TT)^2$,  $(\LL\RR\TT)^2(\LL\RR)^2(\TT\RR)^2$						&$\ZZ_2$\\
$\MN_{17}$ &128&8 &16 &16  & $(\RR\TT\RR\LL)^2$,  $(\LL\RR\TT)^4(\LL\RR)^4$									&$\ZZ_2$\\
$\MN_{18}$ &160&4 &5  &5   & 															&$\ZZ_2^4$\\
$\MN_{19}$ &192&3 &6  &8   & 															&$\ZZ_2^2$\\
$\MN_{20}$ &192&4 &6  &6   & $(\TT(\LL\RR)^2)^3$												&$\ZZ_2^2$\\
$\MN_{21}$ &192&6 &6  &8   &  $\LL\RR\TT\RR\LL\RR\TT\LL\RR\LL(\RR\TT)^2$ 									&$\ZZ_2$\\
$\MN_{22}$ &192&8 &12 &12  & $(\RR\TT)^2(\LL\RR)^2\TT\RR\LL\RR\TT\LL$,  $(\LL\RR)^3\TT\LL(\RR\LL)^2\RR\TT$					&$\ZZ_2$\\
$\MN_{23}$ &192&8 &24 &24  & $(\RR\TT(\RR\LL)^2)^2$,  $(\TT\LL\RR)^3(\LL\RR)^3$									&$\ZZ_2$\\
$\MN_{24}$ &192&8 &24 &24  & $(\RR\TT(\RR\LL)^2)^2$,  $\TT(\LL\RR)^2\TT(\LL\RR)^2\LL\TT\RR\LL\RR$ 						&$\ZZ_2$\\
$\MN_{25}$ &200&4 &4  &10  &  $(\RR\TT\RR\LL)^5$												&$\ZZ_5^2$\\
$\MN_{26}$ &216&4 &6  &12  &   $\TT(\LL\RR\TT\RR)^3$												&$\ZZ_3$\\
$\MN_{27}$ &216&6 &12 &12  &   $\LL(\RR\TT)^2\RR\LL(\RR\TT)^3$, $(\TT\LL\RR\LL\RR)^3$								&$\ZZ_3$\\
$\MN_{28}$ &256&4 &4  &8   & 															&$\ZZ_2$\\
$\MN_{29}$ &256&4 &8  &8   & $(\LL\RR\TT\RR)^2(\LL\RR)^2\LL\TT\RR\LL\RR\TT$									&$\ZZ_2$\\
$\MN_{30}$ &256&4 &16 &16  & $(\RR\TT\RR\LL)^4$, $(\RR\TT\RR\LL(\RR\LL)^2)^2$,$(\LL\RR\TT)^4(\LL\RR)^4$ 					&$\ZZ_2$\\
$\MN_{31}$ &256&4 &32 &32  & $(\LL\RR\TT)^2(\RR\LL)^2(\RR\TT)^2$, $(\LL\RR)^2\TT(\LL\RR)^{14}\TT$						&$\ZZ_2$\\
$\MN_{32}$ &256&4 &64 &64  &  $(\RR\TT\RR\LL)^2$,  $\TT\LL\RR\TT(\LL\RR)^{31}$									&$\ZZ_2$\\
$\MN_{33}$ &256&8 &8  &8   & $(\LL\RR\TT)^2(\LL\RR)^2(\TT\RR)^2$, $\TT(\RR\TT\RR\LL)\TT(\RR\TT\RR\LL)^3$					&$\ZZ_2$\\
$\MN_{34}$ &256&8 &16 &16  & $(\LL\RR\TT\RR\LL\RR\TT)^2,  (\RR\TT)^2\RR\LL(\RR\TT)^2(\RR\LL)^3$							&$\ZZ_2$\\
$\MN_{35}$ &256&8 &16 &16  & $(\LL\RR\TT)^2(\RR\LL)^2(\RR\TT)^2$,$((\RR\TT)^3\RR\LL)^2$,							&$\ZZ_2$\\ 
           &   &  &    &    & $(\LL\RR\TT)^2(\LL\RR)^2\TT(\LL\RR)^3\LL\TT\RR$  &\\
$\MN_{36}$ &256&8 &16 &16  & $((\LL\RR)^2\TT)^2$												&$\ZZ_2$\\
$\MN_{37}$ &256&8 &32 &32  & $(\RR\TT\RR\LL)^2$,   $(\LL\RR\TT)^4(\LL\RR)^{12}$									&$\ZZ_2$\\
$\MN_{38}$ &300&3 &6  &10  & 															&$\ZZ_5^2$\\
$\MN_{39}$ &320&5 &5  &8   & $(\LL\RR\TT\RR)^2\TT(\LL\RR)^2\TT\RR\LL\RR\TT$						 			&$\ZZ_2$\\
$\MN_{40}$ &320&5 &8  &10  &  $(\RR\TT(\RR\LL)^3)^2$, $(\TT\LL\RR)^3\TT\RR(\LL\RR)^2\TT\RR$				 			&$\ZZ_2$\\
$\MN_{41}$ &320&8 &10 &10  & $(\RR\TT)^3(\LL\RR)^4\TT\LL$,$(\TT\LL\RR)^3\LL\RR(\TT\RR)^2\LL\RR$				  			&$\ZZ_2$\\
$\MN_{42}$ &324&3 &6  &18  & $((\LL\RR)^2\TT)^6$   												&$\ZZ_3$\\
$\MN_{43}$ &324&6 &6  &9   &  $(\LL\RR\TT\LL\RR)^2\TT(\LL\RR)^2\TT$, $\TT(\LL\RR(\TT\RR)^2)^3$&$\ZZ_3$\\
$\MN_{44}$ &324&6 &9  &18  &  $(\RR\TT(\RR\LL)^2)^2$, $(\LL\RR\TT)^4\RR\LL(\RR\TT)^2$ 								&$\ZZ_3$\\
$\MN_{45}$ &336&3 &7  &8   &         														&PSL(2,7)\\
$\MN_{46}$ &336&3 &8  &8   & $(\TT\LL\RR)^2(\LL\RR\TT)^2(\LL\RR)^3\LL\TT(\RR\LL)^2\RR$       							&PSL(2,7)\\ 
$\MN_{47}$ &336&4 &6  &8   &   $\TT(\RR\TT\RR\LL)^4$, $(\RR\TT(\RR\LL)^2)^3$ 									&PSL(2,7)\\
$\MN_{48}$ &336&4 &7  &8   &  $(\RR\TT\RR\LL)^3$    												&PSL(2,7)\\
$\MN_{49}$ &336&6 &6  &8   &  $(\LL(\TT\RR)^2)^3$, $(\TT(\LL\RR)^2)^3$										&PSL(2,7)\\
$\MN_{50}$ &336&6 &7  &7   &  $\RR\TT\LL(\RR\TT)^2\RR\LL(\RR\TT)^2$  										&PSL(2,7)\\
$\MN_{51}$ &336&8 &8  &8   & $(\RR\TT\RR\LL)^3$,  $\TT\LL(\RR\TT)^2\LL\RR\TT\RR\LL(\TT\RR)^2$, 							&PSL(2,7)\\
           &   &  &   &    &$(\TT(\LL\RR)^2)^3$&\\ 
$\MN_{52}$ &384&4 &6  &24  & $(\LL\RR\TT)^3(\RR\LL)^2\TT\RR\LL(\RR\TT)^2$									&$\ZZ_2$\\
$\MN_{53}$ &384&4 &12 &24  &  $(\LL\RR\TT\LL\RR)^2\LL\TT\RR\LL\RR\TT$										&$\ZZ_2$\\
$\MN_{54}$ &384&6 &6  &8   &  $(\RR\TT\RR\LL)^3$, $\LL(\RR\TT)^2(\LL\RR)^2\LL(\TT\RR)^2\TT(\LL\RR)^2\TT$  					&$\ZZ_2$\\
$\MN_{55}$ &384&6 &6  &8   &  $(\LL\RR\TT)^3(\RR\TT\RR\LL)^2\RR$ 										&$\ZZ_2$\\
$\MN_{56}$ &384&8 &12 &12  & $(\TT(\LL\RR)^2)^3$,  $((\RR\TT)^3\RR\LL)^2$, &$\ZZ_2$\\
           &   &  &    &    &$(\RR\TT\RR\LL)^4$, $\LL(\RR\TT)^3(\LL\RR)^5\TT$&\\				
$\MN_{57}$ &384&8 &12 &12  &  $((\RR\TT)^3\RR\LL)^2$, $(\RR\TT\RR\LL)^4$, &$\ZZ_2$\\
           &   &  &    &    &$\TT(\LL\RR)^2\TT(\RR\LL)^3\RR\TT\RR\LL\RR$, $\LL(\RR\TT)^3(\LL\RR)^5\TT$&\\ 	
$\MN_{58}$ &384&8 &24 &24  & $\LL(\RR\TT)^2(\LL\RR)^2\TT\RR\LL\RR\TT$, 										&$\ZZ_2$\\
           &   &  &    &    &$((\LL\RR)^3\TT)^2(\LL\RR)^6$&\\  								
$\MN_{59}$ &384&8 &48 &48  & $(\RR\TT(\RR\LL)^2)^2$,  $(\LL\RR\TT)^2\RR\TT\LL\RR\LL\TT(\RR\TT)^2$, 						&$\ZZ_2$\\
           &   &  &    &    &$(\TT\LL\RR)^3(\LL\RR)^9$&\\
$\MN_{60}$ &392&4 &4  &14  &  $(\RR\TT\RR\LL)^7$  												&$\ZZ_7^2$
\end{longtable}
\end{center}

For the maps $\MN_1$ to $\MN_{10}$ detailed descriptions are given in Table \ref{special}.

A {\em genus symbol} is a $6$-tuple $[a,b,c,d,e,f]$ contaning genera of maps $M$, $\Du(M)$,
$\Pe(M)$, $\Pe(\Du(M))$, $\Du(\Pe(M))$ and $\Du(\Pe(\Du(M)))$. If an entry $x$ of a genus symbol is positive, 
then the corresponding map is orientable and its orientable genus is $x$. 
If an entry $x$ is negative then the corresponding 
map is non-orientable and its non-orientable genus is $-x$. An {\em isomorphism symbol} is a 6-tuple $[[a,b,c,d,e,f]]$
that determines which among the maps from the sequence defined above are isomorphic. 
If two entries corresponding to two maps are equal then those maps are isomorphic. 
The {\em hexagonal number} is the number of different entries in an isomorphism symbol.

\begin{center}
\begin{longtable}{|l|c|c|c|l|}
\caption[]{Parallel-product indecomposable non-degenerate reflexible maps $\MN_1$ to $\MN_{10}$ in detail. If 
an underlying graph $G$ has an edge multiplicity $k > 1$, the graph is denoted as $G(k)$.}\label{special} \\

\hline \multicolumn{1}{|l|}{\textbf{Name}} & 
\multicolumn{1}{c|}{\textbf{\textbf{Genus symb.}}} & 
\multicolumn{1}{c|}{\textbf{\textbf{Hex. n.}}} &
\multicolumn{1}{c|}{\textbf{\textbf{Iso. symb.}}} &
\multicolumn{1}{c|}{\textbf{\textbf{Graph}}} 
\\ \hline 
\endfirsthead

\hline \multicolumn{1}{|l|}{\textbf{Name}} & 
\multicolumn{1}{c|}{\textbf{\textbf{Genus symb.}}} & 
\multicolumn{1}{c|}{\textbf{\textbf{Hex.n.}}} &
\multicolumn{1}{c|}{\textbf{\textbf{Iso. symb.}}} &
\multicolumn{1}{c|}{\textbf{\textbf{Graph}}} 
\\ \hline 
\endhead

\hline \multicolumn{5}{|r|}{{Continued on next page}} \\ \hline
\endfoot

\hline \hline
\endlastfoot
$\MN_{1}$  	&  $[0,0,-1,-1,-1,-1]$    	&   3               	&  $[[ 1,1,3,3,5,5]]$ & $K_4$ \\
$\MN_{2}$       &  $[2,2,2,3,2,3]$	  	&   3			&  $[[ 1,2,1,4,2,4]]$ & $C_4(2)$ \\
$\MN_{3}$       &  $[-1,-1,-1,-5,-1,-5]$  	&   3			&  $[[ 1,2,1,4,2,4]]$ & Petersen\\
$\MN_{4}$       &  $[1,1,1,1,1,1]$  		&   1			&  $[[ 1,1,1,1,1,1]]$ & $K_{4,4}$\\
$\MN_{5}$       &  $[3,3,3,5,3,5]$  		&   3			&  $[[ 1,2,1,4,2,4]]$ & $K_{4,4}$\\
$\MN_{6}$       &  $[4,4,4,7,4,7]$  		&   3			&  $[[ 1,2,1,4,2,4]]$ & $C_8(2)$ \\
$\MN_{7}$       &  $[1,1,-5,-5,-5,-5]$  	&   3			&  $[[ 1,1,3,3,5,5]]$ & $DK_{3,3,3}$\\
$\MN_{8}$       &  $[2,2,3,-16, 3, -16]$  	&   6			&  $[[ 1,2,3,4,5,6]]$ & Gen. Petersen $G(8,3)$\\
$\MN_{9}$       &  $[6,6,7,-16, 7, -16]$  	&   6			&  $[[ 1,2,3,4,5,6]]$ & $Q_{3}(2)$\\
$\MN_{10}$      &  $[1,1,1,-11, 1, -11]$  	&   3			&  $[[ 1,2,1,4,2,4]]$ & Pappus
\end{longtable}
\end{center}

\section{Edge-transitive maps}

\begin{figure}[htp]
\centering
\begin{minipage}{0.7\textwidth}
\centerline{\epsfig{file=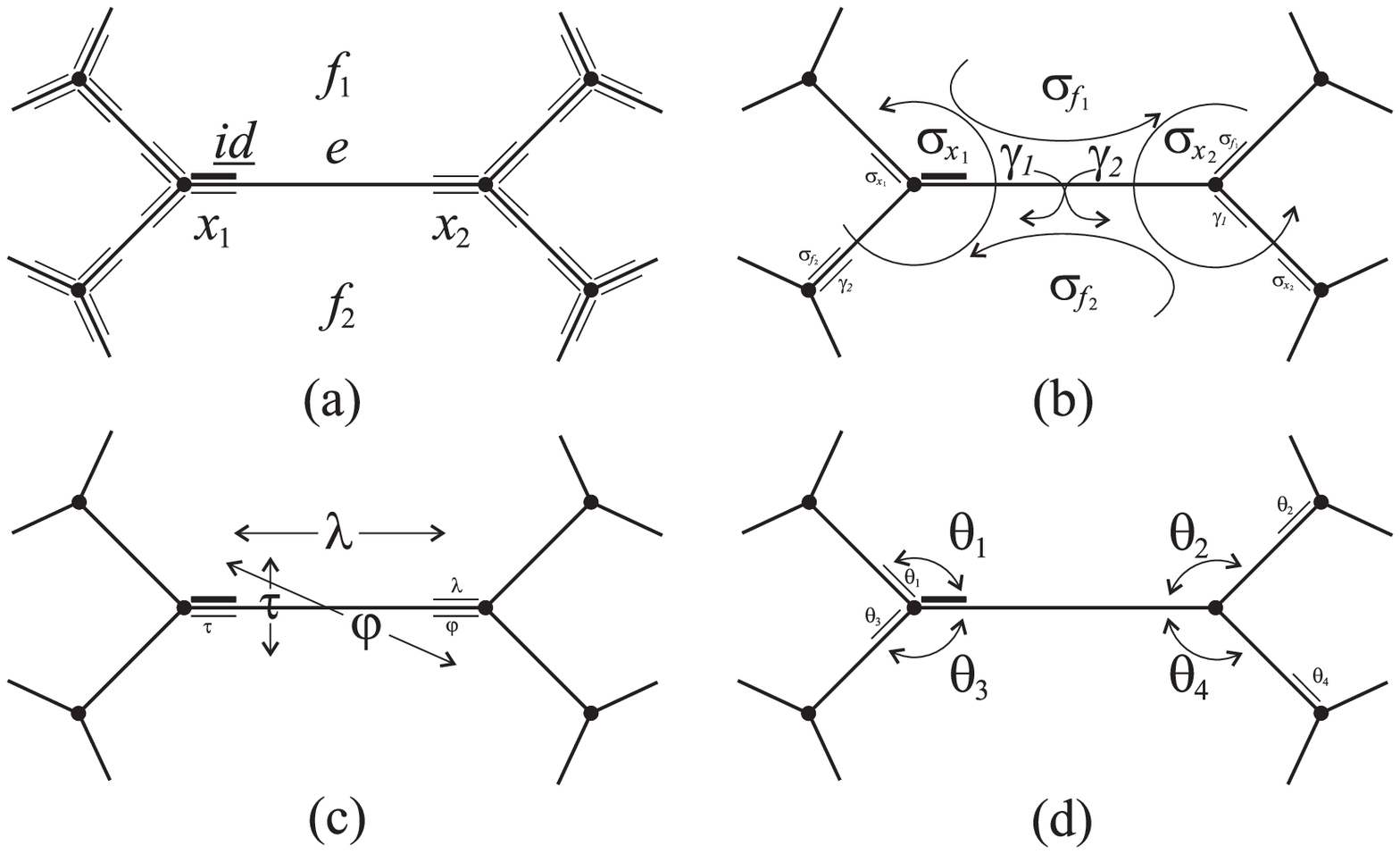,width=\textwidth}}
\end{minipage}%
\hfill
\begin{minipage}{0.3\textwidth}
{\small
\begin{tabular}{|l|l|}
\hline
{\bf Definition} \\
\hline
$\sigma_{x_1} = \alpha_{rt}$    \\
$\sigma_{x_2} = \alpha_{ltrl}$\\
$\sigma_{f_1} = \alpha_{lr}$    \\
$\sigma_{f_2} = \alpha_{trtl}$\\
$\gamma_1 = \alpha_{ltr}$\\
$\gamma_2 = \alpha_{trl}$\\
$\theta_1 = \alpha_{r}$\\  
$\theta_2 = \alpha_{lrl}$\\
$\theta_3 = \alpha_{trt}$\\
$\theta_4 = \alpha_{ltrtl}$\\
$\tau = \alpha_{t}$\\
$\lambda = \alpha_{l}$\\           
$\varphi = \alpha_{lt}$\\
\hline
\end{tabular}
}
\end{minipage}\\[3mm]
\caption{\label{fig_automorphs} Automorphisms "around" the edge $e$ with the root flag $\id$.}
\end{figure}

\begin{figure}[htp]
\centering
\centerline{\epsfig{file=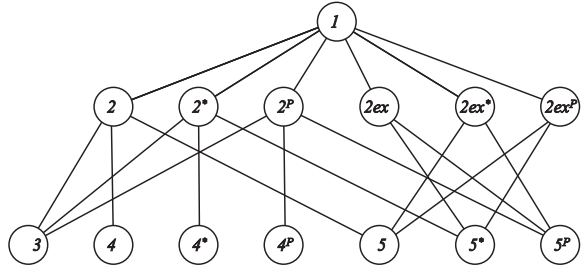,width=0.6\textwidth}}
\caption{\label{fig_order} The partial order of the types of edge-transitive maps. }
\end{figure}

Automorphisms of edge-transitive maps can be studied by focusing on
the situation around the  edge with the root flag. 
In Figure \ref{fig_automorphs}, a set of automorphisms is defined according to how they map
the root flag. In an edge-transitive map not all of those automorphisms are necessarily present.
Let $A$ be the set of all the named automorphisms in Figure \ref{fig_automorphs}. 
Note that those ''named automorphisms'' are not the real automorphisms, but more like the rules how the corresponding 
automorphisms should act, if they exist in an actual map.
For a map $M$, let $A_M$ be a set of all automorphisms from $A$ contained in $M$. Actually, here we have in mind
the set of the  
corresponding automorphisms of the map matching the rules defined by ''named automorphisms'' in $A$. 
According to \cite{memoirs, realizing}, each edge-transitive map can be simply re-rooted, such that
$A_M$ is one of the fourteen {\em edge-transitive types} given in Table \ref{exist_automorphs}. 
Let $A_T$ be a set of automorphisms (''rules'') that a type $T$ map should contain according to Table \ref{exist_automorphs}. 
There is a partial ordering relation $\preceq$ on the set of the types defined by 
$T \preceq T' \Leftrightarrow A_T \subseteq A_{T'}$. 
A Hasse diagram for this ordering is shown in Figure \ref{fig_order}.
The rooting of an edge-transitive map in which the type can be read using Table \ref{exist_automorphs}
is called a {\em canonical rooting}.

An edge-transitive map can have at most two orbits of vertices, faces and Petrie-circuits. The degrees of the vertices in each
of the orbits are denoted by $a_1, a_2$, the sizes of the faces by $b_1$, $b_2$, and the sizes of the Petrie circuits by 
$c_1$, $c_2$. By $\langle a_1, a_2; b_1, b_2; c_1, c_2\rangle$ we denote the {\em map symbol}. If a map is vertex transitive 
then $a_1 = a_2 = a$ and we reduce the symbol to $\langle a; b_1, b_2; c_1, c_2\rangle$. A similar rule extends to faces and
Petrie circuits.

\noindent
\begin{table}[ht]
\centering
\bigskip
\begin{tabular}{|l|l|l|l|}
\hline
{\bf $\Type(M)$} & {\bf $A_M$} & {\bf Map symbol} & {\bf Comments}\\
\hline
1 &  $\tau$, $\lambda$, $\varphi$, $\sigma_{x_1}$, $\sigma_{x_2}$, $\sigma_{f_1}$, $\sigma_{f_2}$, & $\langle a; b; c\rangle$ &\\
& $\gamma_1$, $\gamma_2$, $\theta_1$, $\theta_2$, $\theta_3$, $\theta_4$ & &\\
2 & $\tau$, $\sigma_{x_1}$, $\sigma_{x_2}$, $\theta_1$, $\theta_2$, $\theta_3$, $\theta_4$ &
$\langle a_1, a_2; b; c\rangle$& $2|b, 2|c$\\
$2^*$ &  $\lambda$, $\sigma_{f_1}$, $\sigma_{f_2}$,
$\theta_1$, $\theta_2$, $\theta_3$, $\theta_4$ & $\langle a; b_1, b_2; c\rangle$ & $2|a, 2|c$\\
$2^P$ &  $\varphi$, $\gamma_1$, $\gamma_2$, $\theta_1$, $\theta_2$, $\theta_3$, $\theta_4$ &
$\langle a; b; c_1, c_2\rangle$ & $2|a, 2|b$\\
2ex &  $\tau$, $\sigma_{f_1}$, $\sigma_{f_2}$, $\gamma_1$, $\gamma_2$
& $\langle a; b; c\rangle$ & $a|2$\\
$2^*$ex &  $\lambda$, $\sigma_{x_1}$, $\sigma_{x_2}$, $\gamma_1$, $\gamma_2$,
& $\langle a; b; c\rangle$ & $b|2$\\
$2^P$ex &  $\varphi$, $\sigma_{x_1}$, $\sigma_{x_2}$, $\sigma_{f_1}$, $\sigma_{f_2}$
& $\langle a; b; c\rangle$ & $c|2$\\
3 &  $\theta_1$, $\theta_2$, $\theta_3$, $\theta_4$& $\langle a_1, a_2; b_1, b_2; c_1, c_2\rangle$ & all even\\
4 &  $\sigma_{x_1}$, $\theta_2$, $\theta_4$& $\langle a_1, a_2; b; c\rangle$ & $2|a_1$, $2|a_2$, $4|b$, $4|c$ \\
$4^{*}$ &  $\sigma_{f_1}$, $\theta_3$, $\theta_4$& $\langle a; b_1, b_2; c\rangle$ & $4|a$, $2|b_1$, $2|b_2$, $4|c$\\
$4^{P}$ &  $\gamma_1$, $\theta_2$, $\theta_3$ & $\langle a; b; c_1, c_2\rangle$ & $4|a$, $4|b$, $2|c_1$, $2|c_2$\\
5 &  $\sigma_{x_1}$, $\sigma_{x_2}$& $\langle a_1, a_2; b; c\rangle$ & $2|b$, $2|c$\\
$5^*$ &  $\sigma_{f_1}$, $\sigma_{f_2}$ &$\langle a; b_1, b_2; c\rangle$ & $2|a$, $2|c$\\
$5^{P}$ &  $\gamma_1$, $\gamma_2$& $\langle a; b; c_1, c_2\rangle$ & $2|a$, $2|b$\\
\hline
\end{tabular}
\caption{\label{exist_automorphs}A classification of edge-transitive maps
on 14 types according to the possession of automorphisms around the root flag $\id$.
}
\end{table}

Let us consider a few properties of edge-transitive maps.

\begin{corollary} 
Let $M$ be an edge-transitive map. The product of all 4 simply re-rooted maps is a reflexible map $N$ with
$\Mon(N) = \Mon(M)$ and thus the smallest reflexible cover. 
\end{corollary}
\begin{proof}
Note that from Proposition \ref{totalInfo} it follows that $N \simeq M^M$. 
\end{proof}

The following corollay follows immediately from Proposition \ref{monProject}.
\begin{corollary}
\label{edgeTranProj}
A monodromy quotient  of an edge-transitive map of type $T$ is of type $T'$, such that
$T \preceq T'$.\qed
\end{corollary}

The obvious corollary of Proposition \ref{liftParallel} is the following. 
\begin{corollary}
A parallel product of two canonically rooted edge-transitive maps of type $T$
is an edge-transitive map of type $T'$, such that
$T \preceq T'$. \qed
\end{corollary}

\begin{example}
If we make a parallel product of two edge-transitive maps,
the result need not be edge-transitive. By Proposition \ref{liftParallel}, 
only the lifts of common automorphisms are guaranteed. 
Even, for instance, if we make a parallel product of two simply re-rooted maps 
of type 4, where one map is rooted in $\id$ and the other is rooted in $\id\cdot T$ (relatively to the first map), 
the obtained parallel product in general may not be edge-transitive.
\end{example}

By Corollary \ref{edgeTranProj} the following holds.

\begin{corollary}
If an edge-transitive map of type $T$ is parallel-product decomposable then the factors are 
maps of type $T'$, such that 
$T'\succeq T$.\qed
\end{corollary}

Consider now the impact of the operations $\Du$ and $\Pe$ on  edge-transitive maps.
For the purpose of an easier consideration, we denote $1 = 1^* = 1^P$
and $3 = 3^* = 3^P$.

\begin{proposition}
Let $M$ be an edge-transitive map of type $T$. Then $\Du(M)$ and $\Pe(M)$
are also edge-transitive maps. Furthermore, if $T \in    \{1, 2, 2\mathrm{ex}, 3, 4, 5\}$ then (by abusing the notation),
the types of the map convert as follows:
\begin{align*}
\Du(T) &= T^*, & \Du(T^*) &= T, & \Du(T^P) &= T^P, \\
\Pe(T) &= T, & \Pe(T^*) &= T^P, & \Pe(T^P) &= T^*.
\end{align*}
\end{proposition}
\begin{proof}
Performing the operation $\Du$ can be considered as a renaming (permuting) of the elements $\{\TT, \LL, \TT\LL\}$.
The orbits that determine the edges remain unchanged. The automorphism group $\Aut(M)$ changes the role
and becomes exactly $\Aut(\Du(M))$, but the named automorphisms change their names according to the following.
Let $\dd \in \Aut(F)$ be an automorphism that defines the operation $\Du$. 
Then $\alpha_W \in \Aut(M)$, for $W \in F$, acts
like $\alpha_{\dd(W)}$ in $\Aut(\Du(M))$. For example, $\sigma_{f_2} = \alpha_{trlt} \in \Aut(M)$ acts like 
$\alpha_{\dd(trlt)} = \alpha_{lrtl} = \sigma_{x_2}^{-1} \in \Aut(\Du(M))$. 
A reader can easily verify, that for any type $T$,
the set $A_T$ is in a similar way transformed to a set $\{\gamma_{v_1}^{\pm 1}, \ldots, \gamma_{v_k}^{\pm 1}\}$,
where $A_{T'} = \{\gamma_{v_1}, \ldots, \gamma_{v_k}\}$ and $T'$ is exactly the transformation of the type $T$ as 
claimed in the proposition. A proof for the operation $\Pe$ is similar.
\end{proof}

\begin{corollary}
Each edge-transitive map can be obtained from some map of a type 1, 2, 2ex, 3, 4 or 5 by
one of 6 possible compositions of the operations $\Du$ and $\Pe$.\qed
\end{corollary}

As far as we are considering the analysis of edge-transitive maps 
through their automorphism (and also monodromy) groups, 
we can focus on the types 1, 2, 2ex, 3, 4, and 5. From now on we consider those types only.

From the classification in \cite{memoirs, realizing} the partial presentations of automorphism groups of edge-transitive 
maps can be extracted. 
They are
shown in Table \ref{partial_presentations}. Note that the values of map symbols are used in presentations to denote
partial presentations of maps having a prescribed map symbol. 
The relations, that are independent of a specific map symbol and therefore are present in any partial presentation of the 
corresponding type, are underlined in Table \ref{partial_presentations}. 
The generators and those relations alone determine an {\em universal automorphism group} for the corresponding type.
If a finite presentation of a group $G$ matches the partial presentation corresponding to a type $T$ (for some map symbol), 
we say that $G$ is of the type $T$. This means that $G$ is a finite quotient of the corresponding universal automorphism group.

\noindent
\begin{table}[ht]
\centering
\bigskip
\begin{tabular}{|l|l|}
\hline
{\bf Type} & {\bf A partial presentation for a given map symbol.}\\
\hline
1 & $\langle\tau, \lambda, \theta_1~ \vert\ \underline{ \tau^2, \lambda^2, {\theta_1}^2, (\tau\lambda)^2},
(\theta_1\tau)^a, (\lambda\theta_1)^b, (\tau\lambda\theta)^c, \ldots\rangle$\\
2 & $\langle\tau, \theta_1, \theta_2 ~\vert\ \underline{\tau^2, \theta_1^2, \theta_2^2}, (\theta_1\tau)^{a_1},
(\tau\theta_2)^{a_2}, (\theta_2\theta_1)^{\frac{b}{2}}, (\tau\theta_2\tau\theta_1)^{\frac c 2}, \ldots\rangle$\\
$2$ex & $\langle\tau, \sigma_{f_1}~\vert\ \underline{\tau^2},  (\sigma_{f_1}^{-1}\tau\sigma_{f_1}\tau)^{\frac a 2}, \sigma_{f_1}^b, (\tau\sigma_{f_1})^{c}, \ldots\rangle$\\
$3$ & $\langle\theta_1, \theta_2, \theta_3, \theta_4~ \vert\ \underline{\theta_1^2, \theta_2^2, \theta_3^2, \theta_4^2},
(\theta_1\theta_3)^{\frac{a_1}{2}}, (\theta_4\theta_2)^{\frac{a_2}{2}}, (\theta_2\theta_1)^{\frac{b_1}{2}},$\\
& $~~~~(\theta_3\theta_4)^{\frac{b_2}{2}}, (\theta_4\theta_1)^{\frac{c_1}{2}}, (\theta_3\theta_2)^{\frac{c_2}{2}}, \ldots\rangle$\\
$4$&$\langle\sigma_{x_1}, \theta_2, \theta_4~\vert\  \underline{\theta_2^2,\theta_4^2}, \sigma_{x_1}^{a_1},
(\theta_4\theta_2)^{\frac{a_2}{2}}, (\sigma_{x_1}\theta_4\sigma_{x_1}^{-1}\theta_2)^{\frac{b}{4}},
(\sigma_{x_1}^{-1}\theta_{4}\sigma_{x_1}\theta_{2})^{\frac{c}{4}} , \ldots\rangle$\\
$5$ & $\langle\sigma_{x_1},\sigma_{x_2}~\vert\ \sigma_{x_1}^{a_1},\sigma_{x_2}^{a_2},
(\sigma_{x_1}\sigma_{x_2})^{\frac{b}{2}},(\sigma_{x_1}\sigma_{x_2}^{-1})^{\frac{c}{2}}, \ldots\rangle$\\
\hline
\end{tabular}
\caption{\label{partial_presentations} Partial presentations for automorphism groups of types 1, 2, 2ex, 3, 4, 5. 
}
\end{table}

For a type $T$, a map $M$ is $T$-admissible if there is a subgroup $G \leq \Aut(M)$, such that 
$G$ is generated by automorphisms $A_T$ and no $T' \succ T$ exists, such that automorphisms 
in $A_T'\setminus A_T$
are contained in $G$. In this case, $G$ is called an {\em $T$-admissible} subgroup of $\Aut(M)$. 
Note that $G$ is of type $T$.

To ilustrate the situation here is an example. 
\begin{example}
Take an orientable reflexible (type 1) map $M$. Then the orientation 
preserving subgroup $\Aut^+(M)$ contains and is generated by 
$\{\varphi, \sigma_{x_1}, \sigma_{x_2}, \sigma_{f_1}, \sigma_{f_2}\}$.
But this set is exactly $A_{{\rm 2\mathrm{ex}}^P}$. The subgroup generated by the set does not contain any
other named automorphisms of type $T \succ 2$ex. This is true because
any other automorphism in Figure \ref{fig_automorphs} is not orientation preserving and $M$ is orientable.
Thus $M$ is 2ex$^P$-admissible. Note that saying that a map is 2ex$^P$-admissible in general means that 
a map is either of the type 2ex$^P$ (chiral) or the type 1 (reflexible) and orientable. It is also equivalent to saying
that the map is orientably regular.
\end{example}

It is obvious that every $T$-admissible automorphism subgroup of a map $M$ of type $T$ can 
be represented in a presentation matching the 
corresponding partial presentation of Table \ref{partial_presentations}. 
The following proposition is about the construction of a map from a group in such a presentation.

\begin{proposition}
\label{construct}
Any finite finitely presented group $K$ of type $T$ 
yields the unique $T$-admissible map $M$, such that $K$ is congruent to the $T$-admissible subgroup $G \leq \Aut(M)$.
The construction of $M$ is given
in Table \ref{map_constructions}.
\end{proposition}

\noindent
\begin{table}[ht]
\centering

\bigskip
\begin{tabular}{|l|l|l|l|l|}
\hline
{\bf Type} & {\bf Flags} & {\bf $\TT$} & {\bf $\LL$} & {\bf $\RR$}\\
\hline
1       &  $G$                                      &   $g\cdot \TT = g\tau$                 	&  $g\cdot \LL= g\lambda$       	&   $g\cdot\RR = g\theta_1$\\
\hline
2       &  $G\times\mathbb{Z}_2$                    &   $(g, j)\cdot \TT =$                  	&  $(g, j)\cdot \LL =$               	&   $(g,0)\cdot \RR = (g\theta_1, 0)$\\
        &                                           &    $(g\tau, j)$                   	&  $(g,j+1)$                    	&   $(g,1)\cdot\RR = (g\theta_2, 1)$\\
\hline
$2$ex   &  $G\times\mathbb{Z}_2$                    &   $(g, j)\cdot \TT =$                  	&  $(g, j)\cdot\LL = $              	&   $(g,0)\cdot\RR = (g\sigma_{f_1}^{-1}, 1)$\\
        &                                           &   $ (g\tau, j)$                   	&  $(g, j + 1)$                 	&   $(g,1)\cdot\RR = (g\sigma_{f_1}, 0)$\\
\hline
3       &  $G\times\mathbb{Z}_2\times\mathbb{Z}_2$  &   $(g, j, k)\cdot \TT =$               	&  $(g, j, k)\cdot\LL = $            	&   $(g,0,0)\cdot\RR = (g\theta_1, 0,0)$\\
        &                                           &   $ (g, j+1, k)$                  	&  $(g,j, k+1)$                 	&   $(g,0,1)\cdot\RR = (g\theta_2, 0,1)$\\
        &                                           &                                   	&                               	&   $(g,1,0)\cdot\RR = (g\theta_3, 1,0)$\\
        &                                           &                                   	&                               	&   $(g,1,1)\cdot\RR = (g\theta_4, 1,1)$\\
\hline
4       &  $G\times\mathbb{Z}_2\times\mathbb{Z}_2$  &   $(g, j, k)\cdot \TT = $              	&  $(g, j, k)\cdot\LL = $            	&   $(g,0,0)\cdot\RR = (g\sigma_{x_1}, 1,0)$\\
        &                                           &   $(g, j+1, k)$                   	&  $(g,j, k+1)$                 	&   $(g,0,1)\cdot\RR = (g\theta_2, 0,1)$\\
        &                                           &                                   	&                               	&   $(g,1,0)\cdot\RR = (g\sigma_{x_1}^{-1}, 0,1)$\\
        &                                           &                                   	&                               	&   $(g,1,1)\cdot\RR = (g\theta_4, 1,1)$\\
\hline
$5$     &  $G\times\mathbb{Z}_2\times\mathbb{Z}_2$  &   $(g, j, k)\cdot \TT =$               	&  $(g, j, k)\cdot\LL =$             	&   $(g,0,0)\cdot\RR = (g\sigma_{x_1}, 1,0)$\\
        &                                           &   $(g, j+1, k)$                   	&  $(g,j, k+1)$                 	&   $(g,0,1)\cdot\RR = (g\sigma_{x_2}^{-1}, 1,1)$\\
        &                                           &                                   	&                               	&   $(g,1,0)\cdot\RR = (g\sigma_{x_1}^{-1}, 0,0)$\\
        &                                           &                                   	&                               	&   $(g,1,1)\cdot\RR = (g\sigma_{x_2}, 0,1)$\\
\hline
\end{tabular}
\caption{\label{map_constructions} A construction of the corresponding $T$-admissible map from a partially 
presented group $G$ of type $T$.}
\end{table}

\begin{proof}
If such a rooted map $M$ existed then $\Aut(M)$ would give rise to 
the unique labelling of the flags
in the orbit containing the flag $\id = \id_M$ as follows. Let $G \leq \Aut(M)$ be the  
$T$-admissible subgroup congruent to $K$.
The orbits of $G \leq \Aut(M)$ are blocks of imprimitivity for $\Mon(M)$.
Since $G$ is edge-transitive, there can be at most 4 orbits on the flags and a subgroup 
$Q := \langle \epsilon, \TT, \LL, \TT\LL\rangle$ of order (a most) 4 acts on the set of the orbits transitively. 
Since $Q$ is a small group isomorphic to $\ZZ_2 \times \ZZ_2$, one can easily verify that there is 
always a subgroup $S \leq Q$, such that $S$ acts regularly on the set of the orbits. 
Each flag $x$
can be uniquely
labelled by a pair $(\alpha, w)$, $\alpha \in G$ and $w \in S$, such that $x = \alpha(\id)\cdot w$. 
To see that, let $x = \alpha_1(\id)\cdot w_1 = \alpha_2(\id)\cdot w_2$ for some $\alpha_1, \alpha_2 \in G$ and 
$w_1, w_2 \in S$.
This would imply
$\alpha_2^{-1}(\alpha_1(\id))\cdot  w_1 w_2^{-1} = \id$. Since $\id$ is in the same orbit as $\alpha_2^{-1}(\alpha_1(\id))$ 
and $S$ acts regularly
on the orbits, it first follows $w_1 = w_2$ and then by semi-regularty of $G$ it follows $\alpha_1 = \alpha_2$. 
Thus the labelling is unique and any edge-transitive map corresponds to the unique labelling  $G \times S$. 

The unique labelling alone already determines the map, since a label $(\alpha_W, V)$ 
corresponds to the flag  $\id\cdot WV$. From this information it is straightforward to 
calculate the actions of the involutions 
$\TT$, $\LL$ and $\RR$ on the flags 
with the labels of the form $(\Id, w)$, $w \in S$.
Since for $W \in \{\TT, \LL, \RR \}$,  $x = (\alpha, w)$, $\alpha \in \Aut(M)$, 
it follows $x\cdot W = \alpha(\alpha^{-1}(x)\cdot W)$, the map is uniquely determined by the labelling. 

Note that $S$ is determined by the type of the map. If the type is 1,2, 2ex, 3, 4, 5, then, according to
\cite{memoirs, realizing}, the corresponding sets $S$ are:
$\{\epsilon\}$, $\{\epsilon,\LL\}$, $\{\epsilon,\LL\}$, $\{\epsilon,\TT, \LL, \TT\LL\}$,$\{\epsilon,\TT, \LL, \TT\LL\}$, $\{\epsilon,\TT, \LL, \TT\LL\}$, respectively.
From any finitely presented group $G$ corresponding to  a type $T$, the unique labelling 
$G \times S$ and from that an $T$-admissible map are obtained. The construction for the types following the 
above description is presented in Table \ref{map_constructions}. Here $S$ is modelled by a subgroup of $\ZZ_2\times \ZZ_2$.

For $\alpha, \beta \in \Aut(M)$, $w \in S$, $x = (\alpha, w)$,  it follows that
$\beta(x) = \beta((\alpha, w)) = (\beta\circ\alpha, w)$. It is easy to verify that this is an action.
This action of $G$ on the the labels is consistent 
with the action 
of $G$ on the flags. Also every named automorphism of $G$ maps the root flag exactly according to its name.
Using this, a reader can verify that the maps obtained by the construction from Table \ref{map_constructions} are indeed $T$-admissible.
The conclusion of the proposition follows.

\end{proof}

\begin{corollary}
\label{regularMerge}
Let $M$ be a $T$-admissible map and $G \leq \Aut(M)$ the corresponding $T$-admissible subgroup.
Then the flags of the map $M$ can be partitioned into the blocks of imprimitivity of $G$, such that $G$ acts regularly 
on the blocks.
\end{corollary}
\begin{proof}
Let $G \times S$ be the unique labelling from the proof of Proposition \ref{construct}, where each flag $x \in \Flags(M)$
can be uniquely labelled by $x = (\alpha, w)$, where $\alpha \in G$ and $w \in S$. Then 
$B_{\alpha} = \{x = (\alpha, w)~|~ w \in S\}$, $\alpha \in G$ determine the blocks. The action of $G$ on 
the labels is consistent with the action of $G$ on flags and is defined
as $\beta\cdot B_{\alpha} = B_{\beta\circ \alpha}$. Since $B_{\alpha} = B_{\gamma}$ if and only if $\alpha = \gamma$, 
the action is regular.
\end{proof}

A similar approach in construction of maps from groups using finite presentations was used in \cite{realizing}, 
described in terms of an embedding of an associated Cayley graph in an orientable surface.
One of the problems encountered in \cite{realizing} was whether a finitely presented group matching a  
partial presentation for a type $T$ indeed induces a map of exactly the type $T$. They proved that if the group  
fulfills two conditions, it induces an associated Cayley map of an orientable edge-transitive map of exactly type $T$.
The two conditions were essentially one forcing an orientability and  one preventing other automorphisms
in the obtained map that would imply a type $T' \succ T$.
In Proposition \ref{construct}  a generalized construction to obtain both orientable and non-orientable $T$-admissible 
maps from 
a finite finitely presented group of type $T$ is presented. 
Similar condition for limiting the group automorphisms as Condition 3.2 in \cite{realizing} can be developed and by that
extend some theorems from \cite{realizing} on non-orientable maps.

The author of this work used the programs {\sc Lowx}\cite{lowx} and {\sc Magma}\cite{magma}
to calculate all possible presentations of 
automorphism groups of non-degenerate edge-transitive maps of types 1,2, 2ex, 3, 4, 5 up to 100 edges.  
An edge-transitive map is {\em non-degenerate} if and only if all the values in a map symbol are greater or equal to 3.
During the calculation all possible groups matching the partial presentations form Table \ref{partial_presentations} 
had to be calculated for the type 1 up to size 400, 
for the types 2 and 2ex up to size 200 and for the types 3, 4, 5 up to size 100.
From those presentations one can by Proposition \ref{construct} construct all the corresponding $T$-admissible maps. 
All not $\Aut(M)$-admissible maps were filtered out thus keeping the maps that are of the exact type as the presentation we
started with.

The numbers of triality classes and the numbers of the maps obtained from them for edge-transitive types are 
shown in Table \ref{tab_number}. For the type $1$ (reflexible) and the type 2ex$^P$ (chiral) 
the numbers 
match with Wilson's census of rotary maps \cite{census}.

\noindent
\begin{table}[ht]
\centering
\bigskip
\begin{tabular}{|l|c|c|}
\hline
{\bf Type} 	& {\bf Num. trial. class.} 	& {\bf Num. all. maps}\\
\hline
1 		&  277				& 1223\\
\hline
2 		&  3065				& 16044\\
\hline
2ex		&  66 				& 291\\
\hline
3		&  6033				& 30278\\
\hline
4		&  2980				& 11754\\
\hline
5		&  119				& 495\\
\hline
\end{tabular}
\caption{\label{tab_number} Numbers of triality classes of non-degenerated edge-transitive maps and numbers of all maps that 
can be obtained from the classes using the operations $\Du$ and $\Pe$.
}
\end{table}

Parallel-product decomposition can be applied to edge-transitive maps. The major obstacle to get 
a good characterization (like Theorem \ref{reflexibleDecomposition}) for a parallel-product decomposability of an 
edge-transitive map of type $T$ is the non-regular action of the automorphism and 
the monodromy group. The problem can be solved by changing the presentation of the map, thus also changing the monodromy group.

An {\em universal automorphism group} for an edge-transitive map of type $T$ is any group 
$F = \langle \alpha_{w_1}, \ldots, \alpha_{w_n}~|~ W_1 = \ldots = W_k = 1\rangle$,  
where 
$\{\alpha_{w_1}, \ldots, \alpha_{w_n}\} \leq A_T$ is a set of named automorphisms and 
$\{W_1, \ldots, W_k\}$ is a set of relations, such that any automorphism group of any map $M$ of type $T$ 
is congruent to a quotient
of $F$.

By Proposition \ref{construct}, a $T$-admissible map $M$ is already determined by its $T$-admissible subgroup.
Instead of using the construction in Table \ref{map_constructions} one can work with different presentations
of maps, not in terms of flags but in terms of merged flags, i.e. the blocks described in Corollary \ref{regularMerge}.
But the question is, how should one define a new monodromy group, such that the automorphisms in the usual rooted map presentation
would be also the automorphisms in the new presentation?

Consider the following example.
\begin{example}
Let $F = \langle \tau,\lambda,\rho~|~\tau^2 = \lambda^2, = \rho^2 = (\tau\lambda)^2 = 1\rangle$ be an 
universal automorphism group for the type 1. Let $G$ be a quotient of $F$ that represents an automorphism group of a map $M$,
and $f: F \to G$ the corresponding quotient projection. By Corollary \ref{isomorph}, such a map can be
represented as $M = (f, G, G, 1)$, where the projections of the generators $\tau, \lambda, \rho$ in the quotient
are considered as $\TT$,$\LL$,$\RR$, respectively.

Now we illustrate the correspondence of the actions of the automorphism group and  of the monodromy group. 
Let $N$ be any map of the type 1 (reflexible). Then $\tau(\id) = \id\cdot \TT$. Let $x \in \Flags(N)$. 
By regularity there exists 
the unique $\alpha \in \Aut(N)$, such that $x = \alpha(\id)$. 
Therefore $x\cdot \TT = \alpha(\id)\cdot \TT = \alpha(\id\cdot \TT) = \alpha(\tau(\id))$. Thus if we label the flags 
of $N$ by the  automorphisms, the right action of the monodromy group on the labels 
correspond to the action of $\Aut(M)$ from the right, where $\tau$, $\lambda$, $\rho$ act like $\TT$, $\LL$, $\RR$, 
respectively.
\end{example}
The same concept can be used to define monodromy groups on maps with merged flags, such that the $T$-admissible subgroup
of $\Aut(M)$ acts regularly on the set of merged flags.

Note that this approach matches the concept of a {\em reduced regularity} introduced by A. Breda D'Azevedo 
\cite{reducedBreda} on hypermaps.
Using the concept for defining new kinds of monodromy groups opens a new area of objects to be studied.

From now on, let $F_0 = \langle \tau,\lambda,\rho~|~\tau^2 = \lambda^2 = \rho^2 = (\tau\lambda)^2 = 1\rangle$.
Let 
$F = \langle \alpha_{w_1}, \ldots, \alpha_{w_k}~|~ W_1 = \ldots = W_l = 1\rangle$
be a universal automorphism group for a type $T$, where $\alpha_{w_i} \in A_T$, $w_r \in F_0$ and $W_s \in F$ relations.   
Let $G$ be a finite quotient of $F$ and $f: F \to G$ be the 
corresponding quotient projection.
Define a {\em generalized rooted map in the presentation F} as a quadruple
$M = (f, G = \Mon(M), Z = \Flags(M), \id)$, where $G$ acts transitively and faithfully from the right 
on some finite set $Z$ and $\id \in Z$ is
a root flag. The generators of the monodromy group are exactly
$\{f(\alpha_{i})\}_{i=1}^k \cup \{f(\alpha_{i}^{-1})\}_{i=1}^k$,
i.e. the images of the generators and their inverses 
in the presentation of $F$. Note that the monodromy group together with the chosen set of generators and their inverses
determines the combinatorial and algebraic structure of a generalized rooted map. Therefore, a different choice of 
generators 
in general yields a completely different
combinatorial and algebraic structure. This combinatorial structure is modelled as before by a corresponding colored graph $\Co(M)$
that is the {\em action graph} of $\Mon(M)$ determined by the chosen set of generators and their
inverses. For the theory of action graphs see \cite{malnicAction}.
For generalized map presentations, morphisms, automorphisms, parallel product and quotients are defined in the same way as 
at the beginning of the paper. There we derived all the theory for the special presentation
$F = F_0 = \langle \tau, \lambda, \rho~|~\tau^2 = \lambda^2 = \rho^2 = (\tau\lambda)^2 = 1\rangle$, but
instead of the names of the generators $\tau$, $\lambda$, $\rho$, we used the names $t$, $l$, $r$, respectively.
All the 
claims that did not include the structure of $F_0$
thus hold in general presentations of maps. Note that the only claims that actually 
used the structure of $F_0$ were the claims about the operations
$\Du$ and $\Pe$.
To differ this presentation from others, we will say that the map in this presentation is a map in an {\em usual (rooted) map presentation}.

Note that the concept of reduced regularity can be applied to any map $M$ with 
an usual presentation where merging of 
flags yields blocks of imprimitivity of $\Mon(M)$, such that
some subgroup of $\Aut(M)$ acts regulary on the blocks. 
A new presentation may cause a loss of information meaning that there is no unique construction from a new presentation
to the initial usually presented rooted  map. 
In the edge-transitive case, Proposition \ref{construct} guarantees that the obtained reduced presentations are in one-to-one
correspondence with the corresponding usual rooted map presentations. 

Now let us prove a proposition that links monodromy groups and automorphism groups of generalized rooted maps in some 
presentaion $F$. Note that {\em regular generalized map} means that $\Aut(M)$ is regular on flags.

\begin{proposition}
\label{monAutReg}
Let $M$ be a generalized rooted map in presentation $F$. Then $|\Aut(M)| \leq |\Flags(M)| \leq |\Mon(M)|$. There is equality if and
only if the map is regular. In this case $G := \Mon(M) \simeq \Aut(M)$ and the generalized rooted map $M$ 
is isomorphic to $N = (f, G, G, 1)$.  
\end{proposition}
\begin{proof}
Similarly like for an usual map presentation one can easily verify that $\Aut(M)$ acts semiregularly on flags and 
$\Mon(M)$ is transitive, thus $|\Aut(M)| \leq |\Flags(M)| \leq |\Mon(M)|$.

Let $M$ be regular. To prove regularity of $\Mon(M)$ it suffices to prove that 
the stabilizer $\Mon(M)_{\id}$ is trivial. Let $W \in \Mon(M)$ and $\id\cdot W = \id$. Then for any $d \in \Flags(M)$
there exists an automorphism $\alpha_d \in \Aut(M)$, such that $\alpha_d(\id) = d$. Thus
$$
d\cdot W = \alpha_d(\id)\cdot W = \alpha_d(\id\cdot W) = \alpha_d(\id) = d.
$$
Therefore $W$ is contained in all the stabilizers and thus it is an element of the kernel 
of the action of $\Mon(M)$ acting
on $\Flags(M)$. Since the action is faithful, it follows that $W = \epsilon$ and the action of $\Mon(M)$ is regular.

On the other hand, if $\Mon(M)$ is regular, let $d \in \Flags(M)$. There is an unique element
$W_d\in \Mon(M)$, such that $d = \id\cdot W_d$. Define $\alpha_d(x) = x\cdot W_x^{-1}W_d W_x$. By the regularity of 
$\Mon(M)$, the mapping
is well defined. Let $x \in \Flags(M)$ and $W \in \Mon(M)$. Then
\begin{align*}
\alpha_d(x\cdot W)  &=  \alpha_d(\id\cdot W_x W) = (\id\cdot W_x W) (W_x W)^{-1}W_d (W_x W) =\\
&= \id\cdot W_d W_x W = (\id\cdot W_d W_x)\cdot W = \alpha_d(x)\cdot W.
\end{align*}
It is easy to see that $\alpha_d$ is one-to-one and thus onto.
Thus, $\alpha_d \in \Aut(M)$. Since for every $d \in \Mon(M)$ it follows $\alpha_d(\id) = d$, the group 
$\Aut(M)$ is regular. 

The mapping 
$\gamma: \Mon(M) \to \Aut(M)$, $\gamma: W_d \mapsto \alpha_d$ induces an isomorphism. Since $\id\cdot W_d W_e =$
$\alpha_d(\id)\cdot W_e = \alpha_d(\id\cdot W_e) = \alpha_d\circ \alpha_e(\id)$, the rest follows.
\end{proof}

\begin{example}
To see an example of the use of Proposition \ref{monAutReg}, take a rooted map $M$  of the type 
$2$ex$^P$ (chiral) in an usual map presentation. From Table \ref{partial_presentations} we can see that such a map necessarily contains 
automorphisms $\sigma_{x_1}$ and $\varphi$. It is not hard to see that these two automorphisms generate $\Aut(M)$. Let 
$F = \langle \sigma_{x_1}, \varphi~|~ \varphi^2 = 1\rangle$. Then $G := \Aut(M)$ must be a quotient of $F$ with 
the quotient projection $f: F \to \Mon(M)$. Hence $N = (f, G, G, 1)$  corresponds to the map $M$ but in the presentation 
$F$. Since the type of $M$ is 2ex$^P$, the group $\Aut(M)$ is not regular on $\Flags(M)$. But in the 
new presentation $N$, the same automorphisms yield a regular generalized rooted map. 
If we define $R := f(\sigma_{x_1})$ and $L := f(\varphi)$, we obtain a presentation that is often used when considering
orientably regular maps. Note that the same procedure  applies if $M$ is an orientable reflexible map. In this case the
automorphisms in the new presentation are exactly the original automorphisms that preserve an orientation. Note also, that
a corresponding flag graph is an action graph for generators $R$, $R^{-1}$ and $L$. 
This is the so called truncation of a map.
\end{example}

\begin{example}
\label{type2example}
For the type 2 take $F = \langle \tau, \theta_1, \theta_2 ~|~ \tau^2 = \theta_1^2 = \theta_2^2 = 1\rangle$. 
Any finite quotient of $F$ determines a $2$-admissable map as a regular 
generalized rooted map in the presentation $F$. If $f$ is the corresponding quotient projection, then the 
generators of the monodromy group are $f(\tau)$, $f(\theta_1)$ and $f(\theta_2)$. 
But the monodromy group can be viewed as a monodromy group of some regular hypermap. 
Thus the study of $2$-admissible maps is in a way equivalent to the study of regular hypermaps. 
\end{example}

Since for an edge-transitive map the new monodromy group obtained using the concept of a reduced regularity 
is isomorphic (also congruent) 
to the automorphism group
of the map, the final theorem immediately follows.

\begin{theorem}
\label{theFinal}
An edge-transitive map $M$ is parallel-product decomposable if and only if $\Aut(M)$ contains at least two minimal normal
subgroups. \qed
\end{theorem}

\section{Conclusion}
The main results of the paper are a survey and classification of the quotients of rooted maps, 
the decomposition theorem, its application to the classification of 
reflexible maps of at most 100 edges and its extension to edge-transitive maps.
The necessary presentation theory using the concept of a reduced regularity \cite{reducedBreda} of edge-transitive maps is developed. 
The presentation theory can be extended beyond edge-transitive maps to introduce correspondences
between different combinatorial objects of high symmetry.
For instance, Example \ref{type2example} shows that the classification of edge-transitive 2-admissible maps is about 
as hard as the classification of regular hypermaps. 
The study of several different objects of high symmetry (regular)
is similar and depends only on the presentation of a universal automorphism group. For instance, a theory of 
highly symmetric abstract polytopes can be modelled in a similar way. 
Abstract polytopes have been studied extensively \cite{schulteMullen}. Much less is known about the chiral polytopes or 
other highly symmetric polytopes. 

The decomposition theorem can be used with all such objects and the study of these can be reduced to the 
study of monolithic quotients of the corresponding universal automorphism group. 
Thus the importance of monolithic groups as monodromy groups of parallel-product indecomposable
maps is emphasised. 

Since the algorithm for constructing regular elementary abelian covers or regular maps 
is already developed \cite{elementaryAbelian}, 
a next step could be to specialize that algorithm, so that the group of covering transformations would be 
a monolith in the monodromy group of the cover. Adding this operation to the set of operations $\{\Pe, \Du, \parallel\}$ 
would significantly reduce the set of parallel-product indecomposable maps.  
Another next step would be a study of monolithic groups with a non-abelian monolith. This seems to be a hard problem.

Some of the work in the big paper about Cayley maps \cite{cayleyMaps} can also be extended to generalized rooted maps.
The theory of this paper might also be useful in study of Cayley maps.

Similar approaches using a parallel product, a parallel-product decomposition and the introduced presentation theory can 
be used with any object with a semi-regular
action of the automorphism group, where  the object can be uniquely reconstructed from the group.

\section{Acknowledgement}
I would like to express my gratitude to my supervisors Toma\v z Pisanski and Thomas W. Tucker and to Dragan
Maru\v si\v c for their guidance and support.
I would also like to thank Du\v sanka Jane\v zi\v c and the National Chemical Institute in Ljubljana
for letting me use their computer cluster Vrana.
Also, the author acknowledges the extensive use of 
programs {\sc Lowx} \cite{lowx} and {\sc Magma}\cite{magma}.



\end{document}